\documentclass[12pt,oneside,a4paper,leqno]{article}
\usepackage{amsfonts,amssymb}
\usepackage[centertags]{amsmath}
\usepackage[all]{xy}\xyoption{all}
\usepackage{amsthm}
\newcounter{lemma}
\swapnumbers
\theoremstyle{plain}
\newtheorem{lemma}[equation]{Lemma}\numberwithin{lemma}{section}
\newtheorem{theo}[equation]{Theorem}
\newtheorem{propo}[equation]{Proposition}
\newtheorem{coro}[equation]{Corollary}
\theoremstyle{definition}
\newtheorem{defi}[equation]{Definition}
\newtheorem{example}[equation]{Example}
\newtheorem{remark}[equation]{Remark}
\theoremstyle{remark}

\numberwithin{equation}{section}

\makeatletter
\def\tagform@#1{\maketag@@@{\ignorespaces#1\unskip\@@italiccorr}}

\makeatother

\newcommand{\ldiag}[2]{%
\begin{equation}\label{#1}\begin{aligned}\xymatrix{#2}\end{aligned}\end{equation}%
}
\newcommand{\diag}[1]{%
\begin{equation}\begin{aligned}\xymatrix{#1}\end{aligned}\end{equation}%
}
\newcommand{\ndiag}[1]{%
\begin{equation*}\begin{aligned}\xymatrix{#1}\end{aligned}\end{equation*}%
}
\newdir{ >}{!/8pt/@{ }*@{>}}
\newdir{< }{!/-8pt/@{ }*@{<}}
\theoremstyle{plain}
\newtheorem{pbtheo}[equation]{Pull back theorem}
\newtheorem{bundletheo}[equation]{Bundle theorem}
\newcommand{\hh}{\Phi}

\newcommand{\minus}{\smallsetminus\/}

\newcommand{\R}{\mathbb{R}}
\newcommand{\vect}{{\mathbf{Vect}}}

\newcommand{\V}{{\bf\mathfrak{{F}}}}
\newcommand{\W}{{\bf\mathfrak{{G}}}}

\newcommand{\Vtop}{\mbox{$\V$-${{{\mathbf{Top}}}}$}}
\newcommand{\Diag}{{\mathbf{Diag}}}
\newcommand{\Vdiag}{\mbox{$\V^\op$-${{{\mathbf{Diag}}}}$}}

\newcommand{\homotopic}{\sim}
\renewcommand{\top}{\mathbf{Top}}

\newcommand{\op}{\mathrm{op}}

\newcommand{\from}{\colon}
\newcommand{\st}{\ \mathrm{:} \ }

\newcommand{\pr}{\mathrm{pr}}
\newcommand{\Stra}{\mathbf{Stra}}

\newcommand{\Aut}{\mathrm{Aut}}
\newcommand{\ff}{F} % fibre functor
\renewcommand{\gg}{G} % fibre functor
\newcommand{\hm}{\widehat{\hom}} % hom-bifunctor
\begin{document}
\pagenumbering{arabic}

\title{%
Stratified fibre bundles
}

\author{Hans-Joachim Baues
\footnote{Current address:
Max-Planck-Institut f\"ur Mathematik,
Vivatsgasse, 7 - 53111 Bonn -- Germany,
{\em e-mail: baues@mpim-bonn.mpg.de}.
 } 
 and Davide L.~Ferrario
\footnote{Current address:
Dipartimento di Matematica, 
del Politecnico di Milano,
Piazza Leonardo da Vinci, 32,
20133 Milano -- Italy.
{\em email: ferrario@mate.polimi.it}.
}
}

% \date{\today}
\date{}
\maketitle

\begin{abstract}
A stratified bundle is a fibered space in which strata are classical
bundles and in which attachment of strata is controlled
by a structure category $\V$ of fibers.
Well known results on fibre bundles
are shown to be true for stratified bundles; namely the pull back theorem,
the bundle theorem and the principal bundle theorem.

\emph{AMS SC}:
55R55 (Fiberings with singularities);
55R65 (Generalizations of fiber spaces and bundles);
55R70 (Fiberwise topology);
55R10 (Fibre bundles);
18F15 (Abstract manifolds and fibre bundles);
54H15 (Transformation groups and semigroups);
57S05 (Topological properties of groups of 
homeomorphisms or diffeomorphisms).

\end{abstract}

%\subjclass{
% Primary 
%todo
% Secondary 
%todo
%}

%\keywords{%
%}

%======================================================================

\section{Introduction}
\label{section:introduction}
A stratified bundle is a filtered fibered space $X=\{X_i, i\geq 0\}$
for which the complements $X_i\minus X_{i-1}$, termed strata,
are fibre  bundles. Moreover the attachment of strata is controlled
by a structure category $\V$ of fibers.
For example the tangent bundle of a stratified manifold
is a stratified vector bundle, see \cite{sm}.
Moreover for a compact smooth $G$-manifold $M$ the projection $M\to M/G$
to the orbit space is a stratified bundle by results of Davis \cite{davis}.
In this paper we prove three basic properties:
The \emph{pull back theorem} shows that certain pull backs 
of stratified bundles are again stratified
bundles, see \ref{theo:pullback}.
The \emph{bundle theorem} \ref{theo:bundle}
states that stratified bundles are the same 
as classical fibre bundles in case the structure category
is a groupoid. Moreover the \emph{principal bundle theorem}
\ref{theo:principal}
shows that the function space of admissible maps, $X^V$,
is a stratified bundle and the $\V^\op$-diagram $V\mapsto X^V$
plays the role of the associated principal
bundle, see \ref{theo:principalbundle}.
In the proofs of the results we have to consider the intricate 
compatibility of quotient topology, product topology,
and compact-open topology in function spaces.
As application we study in \cite{sm}
algebraic constructions on stratified vector bundles
(like the direct sum, tensor product, exterior product) and we show
that stratified vector bundles lead to a stratified 
$K$-theory generalizing the Atiyah-Hirzebruch $K$-theory.
This was the main motivation for the proof of the basic
properties in this paper.

%Before describing the structure of the paper, we would like to spend
%a few words on the connections with the existing literature.
%In the last years there have been a significant amount
%of work done on stratified spaces theory in a  topological framework 
%(see  among many others
%\cite{quinn2,farjon2,hugwei,kreck}), without mentioning 
%smooth or semialgebraic structures. 
%Our approach tends to 
%be closer to homotopy theory and algebraic topology than with 
%, so that 
%it might end up with a notation not consistent with the existing
%literature.
% connection with the literature

The structure of the paper is the following. 
In section \ref{section:fibrespaces} the category $\Vtop$ of fibre
families is introduced. 
In section 
\ref{section:cwcomplexes} the proper notion of CW-complex in $\Vtop$
is exploited, and this leads to the definition of stratified bundles
in section \ref{section:stratifiedbundles}.
In more standard language of stratification theory, the stratified 
bundles here introduced should be called CW-stratified bundles.
In the same section 
the pull back theorem \ref{theo:pullback} is introduced.
The proof of this theorem, like the proofs of \ref{theo:bundle},
\ref{theo:principal} and \ref{theo:principalbundle}, is forward-referencing:
its proof will be complete only in section \ref{section:pullback}. 
In section \ref{sec:bundle} we introduce the bundle theorem
\ref{theo:bundle}. Its proof will be done in sections \ref{sec:bundle}
and 
\ref{section:pushout}.
The principal bundle theorem 
\ref{theo:principalbundle} (and also theorem
\ref{theo:principal}) are stated and introduced in section
\ref{section:principalbundle}. Its proof will be the content
of sections \ref{section:9} (where we study NNEP pairs)
and  \ref{section:10} (where NKC categories are introduced).
At the end, in the short section \ref{section:NKC},
some examples of NKC categories are given.
Thus, after this introduction in section \ref{section:introduction},
the paper consists mainly of three parts: 
a first part --- sections \ref{section:fibrespaces} and 
\ref{section:cwcomplexes} ---
with some 
 preliminaries;
then  a second part --- sections \ref{section:stratifiedbundles},
\ref{sec:bundle} and  \ref{section:principalbundle} ---
with our definition of (CW)-stratified bundle and 
the main theorems; finally, a third part  --
sections \ref{section:pullback},
\ref{section:pushout},
\ref{section:9}, \ref{section:10} and \ref{section:NKC} --
 with the actual proofs
of the main theorems, together with some remarks and 
further propositions.

\section{Fibre families}
\label{section:fibrespaces}
Let $\V$ be a small category  
together with a functor $\ff\from \V \to  \top$ to the category 
$\top$ of topological spaces. We assume that the functor $\ff$
satisfies the assumption: 
\begin{quote}
$(*)$ For
every object $V$ in $\V$
the space $\ff(V)$ 
is locally compact, second-countable and Hausdorff. 
\end{quote}
Then $\V$ is 
termed 
\emph{structure category}  
and $\ff$ is a \emph{fibre functor} on $\V$. To be precise
a structure category is a pair
$(\V,\ff)$, but in general we omit to write explicitely $\ff$.
If $\ff$ is a faithful functor then $\V$ is a topological enriched category
in which morphism sets have the compact open topology.
In many examples the functor $\ff$ is actually
the inclusion of a subcategory $\V$ of $\top$ 
so that in this case we need not to mention the fibre
functor $\ff$. 

A 
\emph{fibre family}
with fibres in $\V$ (or a
$(\V,\ff)$-family)
is 
a topological space $X$, termed \emph{total space},
together with a map $p_X\from X \to \bar X$,
termed \emph{projection}
to the \emph{base space} $\bar X$, 
and for every $b\in \bar X$ a selected
homeomorphism $\Phi_b\from p_X^{-1}b \approx \ff X_b$
where $X_b$ is an object in $\V$, called \emph{fibre},
depending on $b\in \bar X$.
The homeomorphism $\Phi_b$ is termed \emph{chart} at $b$.
The family $(p_X\from X \to \bar X, X_b, \Phi_b, b\in \bar X)$
is denoted simply by $X$.
Each object $V$ in $\V$ yields the \emph{point family},
or \emph{$\V$-point},
also denoted by $V$ given by the map $p_V\from \ff V\to *$
where $*$ is the singleton space.

Given two $\V$-families $X$ and $Y$ a \emph{$\V$-map} from $X$ to $Y$
is a pair of maps $(f,\bar f)$ such that the following diagram
\ndiag{% 
X \ar[r]^f
\ar[d]_{p_X} & Y \ar[d]^{p_Y} \\
\bar X \ar[r]^{\bar f} & \bar Y,
}%
commutes, and such that for every $b\in \bar X$ 
the composition given by the dotted arrow of the diagram
\ndiag{% 
p_X^{-1}(b)\ar[r]^{f|p^{-1}(b)}
&  p_Y^{-1}(\bar f b) \ar[d]^{\Phi_{\bar f b}} \\
\ff X_b 
\ar[u]^{\Phi_b^{-1}} 
\ar@{.>}[r]^{} & \ff Y_{\bar f b}
}%
is a morphism 
in the image of the functor $\ff$. That is,
there exists a morphism $\phi\from X_b \to Y_{\bar f b}$
in $\V$
such that the dotted arrow is equal to $\ff(\phi)$.
We will often denote $\ff X_b$ by $X_b$
and it will be clear from the context whether $X_b$ 
denotes an object in $\V$ or a space in $\top$ given by
the functor $\ff$.
If a $\V$-map $f=(f,\bar f)$ is a $\V$-isomorphism
then $f$ and $\bar f$ are homeomorphisms but the converse need not be true.

If $X$ is a $\V$-family  and $Z$ is a 
topological space, then $X\times Z$ is a $\V$-family
with projection 
$p_{X\times Z} = p\times 1_Z\colon X\times Z \to \bar X \times Z$.
The fibre over  a point $(b,z)\in \bar X\times Z$ is equal to 
$p_{X}^{-1}(b)\times \{z\}$; using 
the chart $\Phi_b:p_{X}^{-1}(b) \to X_b$
the chart $\Phi_{(b,z)}$ is defined by $(x,z) \mapsto \Phi_b(x) \in X_b$
where of course we set $(X\times Z)_{(b,z)} = X_b$. 
In particular, by taking $Z=I$ the unit interval
we obtain the cylinder object $X\times I$
and therefore the notion of \emph{homotopy}: two $\V$-maps
$f_0,f_1\colon X \to Y$ 
are $\V$-homotopic (in symbols $f_0\homotopic f_1$)
if there is a $\V$-map $F\colon X\times I \to Y$ 
such that $f_0=Fi_0$ and $f_1=Fi_1$. Here $i_0$ and $i_1$ 
are the inclusions $X \to X\times I$ at the levels $0$ and 
$1$ respectively.

Let $\Vtop$ be the category consisting 
of $\V$-families $p_X: X \to \bar X$ in $\top$ and $\V$-maps.
Homotopy of $\V$-maps yields a natural equivalence relation $\sim$
on $\Vtop$ so that the homotopy category
$(\Vtop)/{_\homotopic}$ is defined.

\begin{remark}
Consider a structure category $\V$ with fibre functor $\ff$.
In general it is not assumed that $\ff$ is faithful, but it is easy
to see that the category $\Vtop$
is equivalent to $\V'$-$\top$,
where $\V'=\V/{_\equiv}$ is the category with the same objects
as $\V$ and morphisms given by equivalence classes
of morphisms, with $f\equiv f' \iff \ff(f) = \ff(f')$. 
We call $\V'$ the \emph{faithful image} category of $\ff$.
\end{remark}

\begin{defi}
If $V$ is an object in $\V$ and $\bar X$ is a space in $\top$,
then the projection onto the first 
factor 
$p_1\from X=\bar X\times V \to \bar X$ 
yields the \emph{product family} 
with fibre $V$;
the charts $\Phi_b\from \{b\}\times V \to V=X_b$
are given by projection and $X_b=V$ for all $b\in \bar X$.
If $X$ is a $\V$-family $\V$-isomorphic to a product
family then $X$ is said to be a \emph{trivial $\V$-bundle}.
In general a \emph{$\V$-bundle} is a locally trivial family of 
fibres, i.e.  a family $X$ over $\bar X$ such that 
every $b\in \bar X$ admits a neighborhood $U$ for which $X|U$ is trivial.
Here 
$X|U$
is the \emph{restriction}
of the family $X$
defined by $U\subset \bar X$.
\end{defi}

Given a family $Y$ with projection $p_Y\colon Y \to \bar Y$ and 
a map $\bar f\colon \bar X \to \bar Y$, the pull-back $X=\bar f^*Y$
is the total space of a family of fibres given by the vertical
dotted arrow of the following pull-back diagram.
\diag{%
X =
\bar f^* Y 
\ar@{.>}[r]
\ar@{.>}[d]
&  Y
\ar[d]^{p}
\\
\bar X 
\ar[r]^{\bar f}
& \bar Y.\\
}%
The charts are defined as follows:  For every $b\in \bar X$ 
let $X_b= Y_{\bar fb}$, and $\Phi_b: p_X^{-1}(b) \to X_b$ 
the composition 
$p_X^{-1}(b)\to p_Y^{-1}(\bar fb) \approx Y_{\bar fb}=X_b$
where the map $p_X^{-1}(b)\to p_Y^{-1}(\bar fb)$
is a homeomorphism since $X$ is a pull-back.

A $\V$-map $i\from A \to Y$ is termed a \emph{closed
inclusion} if $\bar i\from  \bar A \to \bar Y$ is an inclusion, 
$\bar i \bar A$ is closed in $\bar Y$ and the following diagram is 
a pull-back:
\diag{%
i^*Y = A \ar[r]^{i} \ar[d]_{p_A} & Y \ar[d]^{p_Y}\\
\bar A \ar[r]^{\bar i} & \bar Y \\
}%
Hence a closed inclusion
$i\from A \to Y$ induces homeomorphisms
on fibres.

The push-out construction can be extended to the category $\Vtop$,
provided the push-out is defined via  a closed inclusion.
\begin{lemma}\label{lemma:pushout}
Given $\V$-families $A$, $X$, $Y$ and $\V$-maps $f\from A \to X$,
$i\from A\to Y$ with $i$ a closed inclusion 
the push-out diagram 
in $\Vtop$ 
\diag{%
A \ar[r]^f \ar@{->}[d]_{i} & X \ar@{->}[d]\\
Y \ar[r] & Z \\
}%
exists 
and 
$X\to Z$ is a closed inclusion. 
\end{lemma}
\begin{proof}
The push-out is obtained by the following commutative diagram,
\ndiag{
A \ar[rrr]^f 
\ar[rd]^{p_A} 
\ar@{ >->}[ddd]_i 
&  &  & X \ar[ddd] 
\ar[ld]^{p_X}
\\
& \bar A  \ar[r]^{\bar f} \ar@{ >->}[d]^{\bar i} & \bar X \ar[d] \\
& \bar Y  \ar[r]  & \bar Z \\
Y \ar[rrr] \ar[ur]^{p_Y} & & & Z \ar@{.>}[ul]^{p_Z} \\ 
}%
in which the spaces $Z$ and $\bar Z$ are push-outs in $\top$ 
of $f$, $i$ 
and $\bar f$, $\bar i$ respectively; the map $p_Z$ exists and is unique
because of the push-out property. 
Moreover, the charts of fibres of $Z$ are given as follows: 
$\bar Z$  is the union of the two subsets $\bar X$ and $\bar Y\minus \bar A$.
For $b=\bar x \in \bar X\subset \bar Z$ 
there is a canonical composition of 
homeomorphisms
$p_Z^{-1}(b) \approx p_X^{-1}(\bar x) \approx X_{\bar x}$
that yields the chart 
\begin{equation*}
\Phi_{b}:p_Z^{-1}(b)  \to X_{\bar x} = Z_b.
\end{equation*}
In fact, $i:A \to Y$ being a closed inclusion, for every 
$\bar a \in \bar A$ the fibre over $\bar a$ in $A$ is canonically
isomorphic to the fibre over $\bar a$ in $Y$
and hence for every $\bar x \in X$ the fibre over $\bar x$ in $X$
is canonically isomorphic to the fibre over $\bar x$ in $Z$.
Next, for $b=\bar y \in
\bar Y \minus \bar A \subset Z$ 
the homeomorphism is 
\begin{equation*}
\Phi_b: p_Z^{-1}(b) \approx p_Y^{-1}(\bar y) \approx Y_{\bar y} = Z_b.
\end{equation*}
The $\V$-map $j: Y \to Z$ induced by the push-out is a closed inclusion,
since 
$\bar j: \bar Y \to \bar Z$ is the inclusion
of a closed subspace of $\bar Z$; furthermore, $j$ is an isomorphism
on fibres and an inclusion, and therefore  $Y = j^* Z$. 
\end{proof}

A point of this approach is to endow $\Vtop$ with the structure
of a cofibration category (or cylinder category), in order 
to have $\V$-complexes and stratifications. The dual
concept, namely of $\V$-fibrations, was studied by May in \cite{may}.

\section{$\V$-complexes}
\label{section:cwcomplexes}

We recall that for an object $V$ in $\V$ the 
family $\ff(V)\to *$ with base space a singleton is termed a
\emph{$\V$-point} and is also denoted by $V$. 
A disjoint union  of $\V$-points
is called a \emph{$\V$-set}. 
This is a $\V$-family for which the base space has
the discrete topology.
Let $D^n$ be the unit disc in $\R^n$ and $S^{n-1}$ 
its boundary with base point $*\in S^{n-1}$. 
The complement $e^n=D^n\minus S^{n-1}$ is the \emph{open cell}
in $D^n$.
An $\V$-cell is a product family $V\times e^n\to e^n$ with $V\in \V$.

We say that a $\V$-family $X$ is obtained from a 
$\V$-family $D$ by \emph{attaching}  $n$-cells if 
a $\V$-set $Z$ together with a  $\V$-map $f$ is given,
such that the following  diagram
\diag{%
Z\times S^{n-1} \ar[r]^f \ar[d] & D \ar[d] \\
Z \times D^{n} \ar[r]^{\hspace{-24pt} \hh} & X=D\cup_f(Z\times D^n) \\
}%
is a push-out in $\Vtop$.
The inclusion $Z\times S^{n-1} \to Z\times D^{n}$
is a closed inclusion, therefore the push-out exists
and the induced map $D \to X$  is a closed inclusion
and $X\minus D = Z\times e^n$ is a union of open $\V$-cells.
If $Z$ is a $\V$-point then we say that $X$ is obtained from $D$
by attaching a $\V$-cell and $\hh$ is the
\emph{characteristic map} of the $\V$-cell.

\begin{defi}
\label{defi:CWcomplex}
A \emph{relative $\V$-complex} $(X,D)$ is a family $X$ and a 
filtration
\begin{equation*}
D=X_{-1} \subset X_0 
\subset X_1 \subset \dots \subset X_n \subset X_{n+1} \subset \dots \subset X
\end{equation*}
of $\V$-families $X_n$, $n\geq -1$, such that 
for every $n\geq 0$ the $\V$-family $X_n$ is obtained from 
$X_{n-1}$ by attaching $n$-cells and
\begin{equation*}
X = \lim_{n\geq 0} X_n.
\end{equation*}
\end{defi}

The spaces $X_n$ are termed \emph{$n$-skeleta} of $(X,D)$.
If $D$ is empty we call $X$ a 
\emph{$\V$-complex}.
Then $X$ is a union of $\V$-cells.
It is not difficult to show that a relative $\V$-complex $(X,D)$
is Hausdorff and normal provided that $D$ is Hausdorff
and normal. In \ref{defi:CWcomplex} the limit  of skeleta
is the limit induced by the standard limit of topological spaces
(tnd therefore the induced map $X \to \bar X$ is continuous).

\begin{example}
\label{ex:orbit}
Let $G$ be a compact Lie group
and let $\V$ be the \emph{category of orbits} of $G$,
that is, 
$\V$ is the subcategory of $\top$ consisting of spaces
$G/H$, where $H$ is a closed subgroup of $G$,
and $G$-equivariant maps $G/H \to G/H'$.
Then each $G$-CW-complex  (see \cite{tomdieck})
is a $\V$-complex.
\end{example}

\section{Stratified bundles}
\label{section:stratifiedbundles}

Let $A$ be a closed subset  of a space $X$. 
We say that $(X,A)$ 
is a \emph{CW-pair} if there exists a homeomorphism
$(X,A) \approx (X',A')$ of pairs where 
$X'$ is a CW-complex
and $A'$ a subcomplex of $X'$.
A \emph{CW-space} is a space homeomorphic to a CW-complex.

\begin{defi}\label{defi:strati}
We call a 
space $X$ a \emph{stratified} space if 
 a filtration 
\begin{equation*}
X_0 \subset X_1 \subset \dots
\subset X_n \subset \dots \subset \lim_{n\to \infty} X_n = X
\end{equation*}
is given and if 
for every $i\geq 1$ there is  given  a
CW-pair $(M_i, A_i)$ 
and a map
$h_i\from A_i \to X_{i-1}$ with the following properties: 
The subspace $X_i$ is obtained by attaching $M_i$ to $X_{i-1}$ 
via the attaching map $h_i$, i.e.
there is a push-out diagram:
\ndiag{%
A_i \ar@{ >->}[r] 
\ar[d]_{h_i}
& M_i \ar[d] \\
X_{i-1} \ar@{ >->}[r] & X_i.\\
}%
\end{defi}

Moreover $X_0$ is a CW-space.
The complements $X_i\minus X_{i-1}$ are termed \emph{strata}, 
 while the filtration $\emptyset \subset X_0 \subset
X_1 \subset \dots \subset X_n\subset \dots$ is called \emph{stratification} of $X$.
The strata $X_i\minus X_{i-1}$
coincide with the complements $M_i\minus A_i$. 
We call the pairs $(M_i, A_i)$ and the space $X_0$ 
the \emph{attached spaces} of $X$.
If 
%all attached spaces $(M_i,A_i)$ are CW-complexes
%with subcomplexes $A_i$ 
%and if 
all attaching maps are cellular then $X$ is a 
CW-complex and $X_i$ is a subcomplex of $X$.
In this case we say that the stratified space $X$
is a \emph{stratified CW-complex}.
In particular CW-complexes with the \emph{skeletal}
filtration are stratified CW-complexes.
We say that a stratified space is \emph{finite}
if the number of non-empty strata is finite. We always 
assume that a stratified space which is \emph{not finite
is a stratified CW-complex}.
The reason of this assumption is that 
theorems \ref{theo:bundle}, \ref{theo:principal} and 
\ref{theo:principalbundle}
will be proved under the assumption that the base spaces are finite 
or arbitrary CW-complexes: with this choice of terminology
the formulation of the hypothesis  happens to be more compact.
The cellular approximation theorem shows that a stratified 
space is homotopy equivalent to a stratified 
CW-complex.
We point out that a stratified space $X$ 
is a Hausdorff and regular space,
see Gray \cite{gray}.

A map
\begin{equation}
f\from X \to X'
\end{equation}
between stratified spaces is a filtration preserving
map $f=\{f_n\}_{n\ge 0}$ 
together with commutative diagrams
\ndiag{%
X_{i-1}
\ar[d]^{f_{i-1}}
&
\ar[l]
A_i
\ar@{ >->}[r]
\ar[d]
& 
M_i
\ar[d]^{g_i}
\\
X'_{i-1}
&
A'_i
\ar[l]
\ar@{ >->}[r]
&
M_i'
\\
}%
such that $g_i\cup f_{i-1} = f_i$ for $i\geq 1$.
A map is termed \emph{stratum-preserving}
if for every $i$  
\[
f_i (X_i \minus X_{i-1} )
\subset X'_i \minus X'_{i-1}.
\]
Let $\Stra$
be the category of 
stratified spaces  
and stratum-preserving maps.

Consider a manifold with boundary $(M,\partial M)$,
a manifold $N$, and a 
submersion
$h\colon \partial M \to N$.
The push-out of $h$ and the inclusion $\partial M \subset M$
\diag{
\partial M \ar@{ >->}[r] 
\ar[d]_h & M 
\ar[d]
\\
N \ar[r] & 
M \cup_h N.\\
}%
yields a stratified space $X=M \cup_h N$ with stratification
$X_0=N \subset X_1=X$.
For example if a diffeomorphism $\partial M \cong Z \times P$ is given
where $Z$ and $P$ are manifolds and if $h:\partial M \cong Z\times P \to Z$
is defined by the projection then $X=M\cup_hN$ is a 
\emph{manifold with singularities}, see 
Rudyak \cite{rudyak}, Baas \cite{baas},
Botvinnik \cite{botvinnik},
Sullivan \cite{sullivan},
Vershinin \cite{vershinin}.
Also the stratified manifolds (stratifolds) of Kreck \cite{kreck}
are stratified spaces.
Manifolds with singularities and stratifolds have ``tangent bundles''
which are stratified vector bundles. These are important examples
of stratified fibre bundles introduced in the next definition.

Recall that $\V$ denotes a structure category together 
with a fibre functor $\ff\from \V\to \top$.
\begin{defi}\label{defi:stratifiedvectorbundle}
A \emph{$\V$-stratified fibre bundle} 
is a stratified space $\bar X$ together with a $\V$-family
\begin{equation*}
X \to \bar X
\end{equation*}
with the following properties. For $i\geq 1$ the restriction 
$X_i = X|\bar X_i = X_{i-1} \cup_{A_i} M_i$ is the push-out of $\V$-maps
\ndiag{%
M_i
\ar[d]
&
\ar@{ >->}[l]
A_i
\ar[d]
\ar[r]^{h_i}
&
X_{i-1} 
\ar[d]
\\
\bar M_i
&
\bar A_i
\ar@{ >->}[l]
\ar[r]^{\bar h_i}
&
\bar X_{i-1} 
\\
}%
where $M_i \to \bar M_i$ is a $\V$-bundle and $A_i=M_i|\bar A_i$.
Moreover, $X_0 \to \bar X_0$ is a $\V$-bundle
and $X=\lim_{i\to \infty} X_i$.
Hence the strata $X_i \minus X_{i-1} \to \bar X_i \minus \bar X_{i-1}$
are $\V$-bundles.
\end{defi}

\begin{remark}
It is possible to show that a finite $\V$-stratified fibre bundle
yields a 
\emph{stratified system of fibrations} (according to the definition 
8.2 at page 420 of 
 F. Quinn in \cite{quinn2}) on $\bar X$.
In fact, since $\bar X$ is a CW-complex and hence its
strata are paracompact, the restrictions of $p_X \from X \to \bar X$
to the strata are  fibrations (theorem XX.4.2.(2) at page 405 of
\cite{dugundji} -- see also 
 theorem 1.3.5 at page 29  and
 Ex. 1 at page 33 of \cite{fripic}).
Furthermore, for every $n\geq 0$ the inclusion $X_n \subset X$
is a $\V$-cofibration and all the spaces are compactly generated, 
therefore it is possible to prove that 
$\bar X_n$ is a $p$-NDR subset of $\bar X$,
in the notation of \cite{quinn2}. Thus a $\V$-striatified 
bundle $X$ yields in a natural way a stratified  system
of fibrations on $\bar X$.
Actually, if the CW-structure of $\bar X$ is given by 
a triangulation or a simplicial structure,
then $X\to \bar X$ yields  what in \cite{quinn2,farjon1,farjon2}
is termed a \emph{simplicially stratified fibration}. 

Starting from the papers  \cite{quinn1,hughes1,hughes2},
F. Quinn and B. Hughes 
developed a theory of homotopically stratified sets
which goes in a different direction than this paper.
In their approach isotopy Whitney lemmas and applications
to surgery theory are central, while our motivations
reside more in abstract homotopy theory  and $K$-theory (see
\cite{sm} and  \cite{baues}). 
For a general and modern introduction 
to stratified spaces and stratified bundles,
related to problems in surgery theory 
a common reference is \cite{hugwei}.
\end{remark}

\begin{example}
By comparing definition \ref{defi:stratifiedvectorbundle} 
and definition 
\ref{defi:CWcomplex}
it is easy to see that 
a $\V$-complex with the skeletal filtration is a
$\V$-stratified bundle since the spaces
$(V\times D^n,V\times S^{n-1})$ are trivial $\V$-bundles
over the CW-pair $(D^n,S^{n-1})$.
Moreover, if $X$ is a $\V$-stratified bundle such that 
all the attaching maps are cellular, 
then $X$ is a $\V$-complex.
Here we use the bundle theorem in section \ref{sec:bundle} below.
In general, a  $\V$-stratified bundle 
has the $\V$-homotopy type of a $\V$-complex
(this is a consequence of a cellular approximation theorem
for $\V$-complexes, see \cite{baues}).
\end{example}

\begin{example}
\label{ex:davis}
Let $G$ be a compact Lie group and $M$ a compact smooth $G$-manifold.
Let $\V$ be the orbit category of $G$.
Then the  \emph{augmented $\mathcal{G}$-normal system}
associated to $M$, 
as defined by Davis in \cite{davis},
yields naturally a $\V$-stratified bundle $M\to \bar M = M/G$
via the corresponding assembling functor.
The strata are defined exactly as in the proof of theorem 4.9
of \cite{davis}, by $M_n = M(n-1)$ and $A_n= \partial M_n$.
Thus each open stratum contains the disjoint union 
of all the open strata given by the stratification
of $M$ (or the corresponding stratification on $\bar M$) 
by normal orbit type of depth $n$ in the poset 
of the normal orbit types.
\end{example}

A \emph{$\V$-stratified map} $f\from X \to X'$ 
between $\V$-stratified bundles is 
is  given by sequences
$\{f_n\}_{n\geq 0}$ and $\{g_n\}_{n\geq 1}$ of $\V$-maps
such that, given the commutative diagrams 
\ndiag{%
X_{i-1}
\ar[d]^{f_{i-1}}
&
\ar[l]
A_i
\ar@{ >->}[r]
\ar[d]
& 
M_i
\ar[d]^{g_i}
\\
X'_{i-1}
&
A'_i
\ar[l]
\ar@{ >->}[r]
&
M_i'
\\
}%
for every $i$, we have 
$g_i\cup f_{i-1} = f_i$ for $i\geq 1$.
If for every $i$ the inclusion 
$f(X_i\minus X_{i-1}) \subset X'_i \minus X_{i-1}'$  
holds, $f$ is termed \emph{stratum-preserving}.

\begin{pbtheo}\label{theo:pullback}
Let $\bar X$ and $\bar X'$ be finite stratified spaces with attached 
spaces which are locally finite and countable CW-complexes.
Let $\bar f\from \bar X \to \bar X'$
be a stratum-preserving map (in $\Stra$) and $X'\to\bar X'$
a $\V$-stratified bundle. 
Then the pull-back $\bar f^*X' \to \bar X$ is a $\V$-stratified bundle.
\end{pbtheo}
\begin{proof}
Let $n$ be the number of strata of $\bar X$
and $X_i = \bar f^* X_i'$ for $i=1,\dots n$. 
Clearly $\bar f_0^* X_0'$ is a $\V$-bundle on $\bar X_0$.
Since $X'$ is a $\V$-stratified bundle 
there are $\V$-bundles $M_i'$ and 
$\V$-maps $h'_i\from A'_i \subset M_i' \to X_{i-1}'$
such that 
$X_{i}' = X_{i-1}'\cup_{h_i'} M_i'$.
Let the attached spaces of $\bar X$ be
denoted by $(\bar M_i, \bar A_i)$ and 
$\bar h_i$ the attaching maps.
The stratified map  $\bar f$ yields maps
$\bar g_i\from \bar M_i \to \bar M_i'$
such that $\bar g_i \cup \bar f_{i-1} = \bar f_i$.
For every $i\geq 1$ let $M_i$ be the  pull-back
$\bar g_i^* M_i' \to \bar M_i$.
It is a $\V$-bundle. Let $A_i$ be its restriction
to $\bar A_i$.
By assumption 
$\bar M_i$ and $\bar M'_i$ are locally finite
and countable; furthermore since $\bar f$ is stratum-preserving,
for every $i$ we have 
$\bar f(\bar X_{i} \minus \bar X_{i-1}) 
\subset \bar X'_i \minus \bar X'_{i-1}$.
Thus, 
by applying proposition \ref{propo:theo:pushout2} below, the 
pull-back family $X_i$
is obtained by attaching the $\V$-bundle
$M_i$ to
$X_i$ via a $\V$-map $h_i\from A_i\subset M_i \to X_i$,
where $h_i$ is a suitable $\V$-map induced by the construction.
This is true for $i=1,\dots, n$, 
hence $X_n = X$ is a $\V$-stratified bundle.
\end{proof}

A similar theorem was proved by Davis \cite{davis} (theorems 1.1 and 1.3),
in the case of the pull-back $f^*M$ of a smooth $G$-manifold  $M$,
where $f$ is a (weakly) stratified map and $G$ is a compact Lie group.

Let $(\V,\ff)$ and $(\V',\ff')$
be structure categories. Then we define the structure 
category $(\V\times\V',\ff\times\ff')$ by the functor
\[
\ff\times\ff'\from \V\times\V' \to \top.
\]
Here $\V\times\V'$ is the product category consisting
of pairs $(V,V')$ of objects $V\in \V$ and $V'\in\V'$
and pairs of morphisms.
The functor $\ff\times\ff'$ carries $(V,V')$
to the product space
$\ff(V)\times\ff'(V')$.
For a $\V$-family $X$ and a $\V'$-family $X'$ with the 
base space $\bar X=B=\bar X'$ we define the \emph{fiberwise product}
$X\times_BX'$ by the pull-back diagram in $\top$
\ndiag{%
X\times_BX' \ar[r]\ar[d] \ar@{.>}[rd]^{p} & X' \ar[d] 
\\
X \ar[r] &
B \\
}%
Then it is clear that $X\times_BX'\to B$ is a 
$(\V\times \V',\ff\times\ff')$-family.
Moreover we get compatibility with stratifications as follows.

\begin{propo}
Let $X$ be a $\V$-stratified bundle and let $X'$ be a 
$\V'$-stratified bundle such that the stratified spaces
$\bar X = B = \bar X'$ coincide.
Moreover assume that $X$ and $X'$  are locally compact.
% Hausdorff is automatical
% by Cohen's theorem (9.4 pag. 248 Dugundji)
% a space is compactly generated if and only if
% it is the quotient of a locally compact space.
% that is, a locally compact (X x X') is always compactly generated 
% (if Hausdorff, of course).
%%%and that $X\times_BX'$ has the compactly generated 
%%%topology.
Then $X\times_BX'\to B$
is a $\V\times \V'$-stratified bundle.
\end{propo}
\begin{proof}
By construction $X$ is the quotient space of a projection
\[
q\from X_0 \cup \coprod_{i\geq 0} M_i \to X,
\]
and the same holds for $X'$,
with a projection $q'$.
By the Whitehead theorem (see \cite{dugundji}, theorem
XII.4.1 at page 262) 
the map $q\times q'$ is a quotient map onto $X\times X'$,
since $X$ and $X'$ are locally compact.
Now, $X \times_B X \subset X\times X'$
is a closed subspace, hence 
the restricted projection
\[
q\times q'\from (q\times q')^{-1} X\times_B X \to X' \times_B X'
\]
is a quotient map.
We can consider now the union
\[
(X_0\times_{\bar X_0} X_0')\coprod_i (M_i\times_{\bar M_i} M_i') 
\subset (q\times q')^{-1} X\times_B X \to X'.
\]
Then it is possible to see that the projection $q\times q'$
yields a quotient map
\[
(X_0\times_{\bar X_0} X_0')\coprod_i (M_i\times_{\bar M_i} M_i' )
\to
X \times_B X',
\]
which gives the $\V\times \V'$-stratification of $X\times_B X'$.
\end{proof}

\section{The bundle theorem}
\label{sec:bundle}
A \emph{groupoid} is a category in which all morphisms are isomorphisms.
If the structure category $\V$ is a groupoid then $\V$-complexes
and $\V$-stratified bundles are actually $\V$-bundles.
More precisely we show the following result which we could not
find in the literature though special cases are well known
like the clutching construction of Atiyah  \cite{at}, page 20.
This can be briefly explained as follows: 
if $\bar X=\bar X_1 \cup \bar X_2$ and $\bar A = \bar X_1 \cap \bar X_2$,
all the spaces are compact and $E_i \to \bar X_i$ 
are (vector) bundles with an isomorphism 
of the restrictions $\varphi\from E_1|A  \cong E_2|A$, 
then the clutching (or glueing) construction is a way
to define the push-out of the two bundles $E_1\cup_\varphi E_2$
over the push-out space $\bar X_1 \cup_{\bar \varphi} \bar X_2$.

\begin{bundletheo}
\label{theo:bundle}
Let $\V$ be a structure category which is a groupoid.
Then a $\V$-stratified bundle $X \to \bar X$ 
is a $\V$-bundle over $\bar X$. Conversely,
given a stratified space $\bar X$ and a $\V$-bundle
$X$ over $\bar X$ then $X$ is a $\V$-stratified bundle.
In particular, a $\V$-complex $X\to\bar X$ 
is a $\V$-bundle over $\bar X$ and,
given a CW-complex $\bar X$ and a $\V$-bundle
$X$ over $\bar X$, then $X$ is $\V$-isomorphic to a  $\V$-complex.
\end{bundletheo}
\begin{proof}
If $X$ is a $\V$-complex then
the theorem is a consequence of corollary 
\ref{coro:clutch2} and lemma \ref{lemma:cwpair} below.
Otherwise, assume that $X$ is a $\V$-stratified bundle.
Since $X$ is obtained by attaching $\V$-bundles
to $\V$-bundles via $\V$-maps,
we can apply a finite number of times 
corollary \ref{coro:clutch2} below and obtain the result.
On the other hand 
assume that $\bar X$ is a stratified space and that $X\to \bar X$ 
is a $\V$-bundle.
By applying a finite number of times 
lemma \ref{lemma:attaching} it is possible to show 
that $X$ is a $\V$-stratified bundle.
\end{proof}

\begin{example}\label{ex:3}
Let $\V$ be the category consisting of finite sets $\{1,2,\dots,n\}$
and permutations of these sets.
A $\V$-family $X$ is a finite-to-one map $X\to \bar X$.
By theorem \ref{theo:bundle}, 
if $X$ is a $\V$-complex, then 
$p_X$ is locally trivial and  hence a covering map.
\end{example}

%(Compare with the \emph{clutching} construction, in Atiyah \cite{at}
% page 20).

\begin{lemma}\label{lemma:clutch1}
Let $\V$ be a groupoid.
Let $Y\to \bar Y$ be a $\V$-bundle, $\bar Z$  a set
and $h\from V\times \bar Z \times S^{n-1} \to  Y$ 
a $\V$-map  which is the attaching map 
of $X=Y\cup_h (V\times \bar Z \times D^n)$.
Given a set $\bar U \subset \bar Y$ open in $\bar Y$
such that $Y|\bar U$ is trivial,
there is an open set 
$\bar U'' \subset \bar X$ such that $\bar U'' \cap \bar Y = \bar U$
and such that $X|\bar U''$ is trivial.
As a consequence, the push-out space $X=Y\cup_h (V\times \bar Z\times D^n)$
is a $\V$-bundle.
\end{lemma}
\begin{proof}
Let $\hh$ denote the characteristic 
map $\hh\from V\times \bar Z\times D^n \to X$
and $\bar \hh$, $\bar h$ the maps induced on the base spaces.
The restriction of $V\times \bar Z \times S^{n-1}$ 
to $\bar h^{-1}\bar U$ is trivial,
and with a suitable change of coordinates
(which exists since $\Aut_\V(V)$ is a topological
group)
we can
assume that the restriction
\[
h_0 \from V \times \bar h^{-1}\bar U \to Y|\bar U \approx V\times \bar U
\]
is of the form $1_V\times \bar h$.
Now, there is an open set $\bar U' \subset  \bar Z\times D^n$
with the property that 
$\bar U' \cap  \bar Z\times S^{n-1} = \bar h^{-1}\bar U$.
Since $\bar X = \bar Y \cup_{\bar h} (\bar Z \times D^n)$
the image $\bar U'' = \bar \hh (\bar U') \cup \bar U$ is open in $\bar X$.
Moreover, since $V$ is locally compact Hausdorff 
the following diagram is a push-out.
\ndiag{%
V \times \bar h^{-1}\bar U 
\ar[r]
\ar[d]
&
V \times \bar U \ar[d] 
\\
V\times \bar U' 
\ar[r]
&
V \times \bar U''\\
}%
But by assumption 
$X$ is the push-out $Y\cup_h(V\times \bar Z \times D^n)$,
and the restriction of the push-out to the pre-images 
of $\bar U''$ yields a push-out
\ndiag{%
V \times \bar h^{-1}\bar U 
\ar[r]
\ar[d]
&
V \times \bar U 
\ar[d] 
\\
V\times \bar U' 
\ar[r]
&
X|\bar U''.\\
}%
Therefore $X|\bar U'' \approx V \times \bar U''$
by an isomorphism which is an extension of the chosen
isomorphism $Y|\bar U \approx V \times \bar U$.
This implies that every point in $\bar Y$ has a neighborhood 
$\bar U''$ in $\bar X$ such that  $X|\bar U''$ is trivial.
On the other hand, if $x$ is a point in $\bar X\minus \bar Y$,
then there is a neighborhood $\bar U''$ of $x$ contained 
in $\bar X \minus \bar Y$
and in this case $X|\bar U''$ is trivial
because it is $\V$-homeomorphic 
to $ V\times \bar \hh^{-1} \bar U'' \subset  V\times \bar Z \times e^n$.
\end{proof}

\begin{coro} \label{coro:clutch2}
Let $(X,Y)$ a relative $\V$-complex.
Assume that $\V$ is a groupoid 
and that $Y \to \bar Y$ is a $\V$-bundle.
Then $X\to \bar X$ is a $\V$-bundle.
\end{coro}
\begin{proof}
Consider the cellular filtration 
\[
Y=X_{-1} \subset X_0 \subset X_1 \subset  \dots \subset X_n 
\subset \dots \subset X=\lim_{n} X_n.
\]
By applying inductively lemma \ref{lemma:clutch1}
we see that $X_n$ is a $\V$-bundle over 
$\bar X_n$ for every $n\geq 0$.
Moreover, there is a sequence of 
sets $\bar U_n \subset \bar X_n$ 
open in $\bar X_n$ such that 
$\bar U_n \cap \bar X_{n-1} = \bar U_{n-1}$ 
and  $\V$-isomorphisms 
$\alpha_n\from X_n | \bar U_n \approx V \times \bar U_n$.
Moreover $\alpha_n$ can be chosen to be extension of $\alpha_{n-1}$
(see the proof of \ref{lemma:clutch1}).
Since $V$ is locally compact Hausdorff we can 
take limits, 
and obtain a $\V$-homeomorphism
\[
X|\bar U \approx V \times \bar U,
\]
where $\bar U = \lim_{n} \bar U_n$.
\end{proof}

\section{The principal bundle theorem}
\label{section:principalbundle}

It is well known that each fibre bundle $X\to \bar X$
with fibre $V$ and structure group $G$ yields the associated
principal bundle $X_G\to \bar X$ with fibre $G$ such that 
$X_G$ is a right $G$-space for which there is an isomorphism
of bundles.
\ldiag{diag:5:1}{%
X \ar[rr]^{\cong} \ar[dr] & & X_G\times_G V \ar[dl]\\
& \bar X\\
}%
Moreover Steenrod \cite{steenrod} p.~39 points out that in case $G$ 
is a subgroup of the group of homeomorphisms of $V$ with the compact
open topology then the principal bundle 
\begin{equation}\label{eq:steenrod}
X_G = X^V
\end{equation}
is the function space $X^V$ of ``admissible'' maps $V\to X$
with the compact open topology.
The action of $g\in G$ on $\alpha\in X^V$ is given by
composition
\[
\alpha \cdot g = \xymatrix@1{V\ar[r]^g & V \ar[r]^\alpha & X \\} \in X^V.
\]
In this section we generalize these results for $\V$-complexes
and $\V$-stratified bundles.

\begin{defi}\label{defi:phiV}
Given a faithful fibre functor $\ff\from \V \to \top$ and an object $V\in \V$
we obtain a new fibre functor
\[
\ff^V\from \V \to \top
\]
which carries an object $W$ in $\V$ to the space $\ff^V(W) = \hom_\V(V,W)$
with the compact open topology.
The functor $\ff^V$ need not be faithful, even though $\ff$ is faithful.
\end{defi}

\begin{remark}
Since for every $V$ the fibre $\ff(V)$ is 
Hausdorff, we see that for each $W\in \V$ the fibre $\ff^V(W)=W^V$ 
of $\ff^V$ is Hausdorff.
Moreover, since  $W$ and $V$ are $2^\circ$ countable, also
$W^V$ is $2^\circ$ countable (see e.g. \cite{dugundji}). 
We will show that 
if $(\V,\ff)$ is a NKC structure category
%(see definition\ref{defi:NKC} below)
then $W^V$ is locally compact so that $\ff^V$ fulfills
assumption $(*)$.
\end{remark}

The definition of a NKC structure category is the following.
First, 
a family $\mathcal{K}$ of compact sets of $V$ is termed
\emph{generating} if for every 
$\V$-family $Y$ 
the subsets 
\[
N_{K,U} = \{ f\in Y^V, f(K)\subset U\}
\]
with $K\in \mathcal{K}$
and $U$ open in $Y$ yield a sub-basis for the topology of $Y^V$.

\begin{defi}\label{defi:NKC}
A structure category $\V$ with faithful fibre 
functor $\ff$ has the \emph{\emph{NKC}-property}
if 
for every $V\in \V$
there is a generating family of compact sets
$\mathcal{K}$
such that for every $K\in \mathcal{K}$, every 
$W\in \V$ and every compact subset $C\subset W$ 
the subspace
\[
N_{K,C} \subset \hom_\V(V,W) = W^V
\]
is compact.
Examples are given in section \ref{section:NKC} below.
\end{defi}

\begin{remark}
\label{rem:NKC}
Since the spaces 
$W$ in $\V$ 
are assumed locally compact
this implies that for  a NKC category 
the function spaces $W^V = \hom_\V(V,W)$ are locally compact.
Furthermore, it is not difficult to see that a closed 
subcategory of a NKC category is also a NKC subcategory.
\end{remark}

\begin{example}\label{ex:1}
Let $\V$ denote the category consisting of 
a finite dimensional vector space $V$ and morphisms given by
elements in a closed subgroup $G$ of $GL(V)$. The fibre functor 
is the embedding $\ff\from \V \to \top$.
Then,
by \ref{theo:bundle}, 
a $\V$-complex $X$ is a vector bundle with structure group $G$
and fibre $V$. We will see in  theorem \ref{theo:NKC}
that $(\V,\ff)$ has the NKC property.
Therefore by theorem \ref{theo:principal}
the principal bundle $X^V$ is a $(\V,\ff^V)$-complex 
and hence a $(\V,\ff)$-bundle by \ref{theo:bundle}.
This is the result of Steenrod in \ref{eq:steenrod}.
\end{example}

\begin{example}\label{ex:2}
Let $G$ be a compact Lie group and $\V$ its orbit category
(see \ref{ex:orbit}). The fibre functor
$\ff\from \V \to \top$ is the embedding.
Then a smooth $G$-manifold $X$ has a structure of a $\V$-complex
where $\bar X=X/G$ and $p_X \from X \to \bar X$ is the projection
onto the orbit space.
For every closed subgroup $H\subset G$ the orbit $V=G/H$
is an object of $\V$ and the function space 
$X^V$ is homeomorphic to the fixed subspace $X^H \subset X$
via the evaluation map at $1\in G$. 
We will see as a consequence of proposition \ref{propo:theo:allcompact} 
that the orbit category $\V$ with the embedding functor
$\ff\from \V \to \top$ has the NKC property. Hence 
by theorem \ref{theo:principal} the space $X^H$ is a $\V$-complex
with fibre functor $\ff^{G/H}$.
\end{example}

We now define for every $(\V,\ff)$-family $Y$ the
associated $(\V,\ff^V)$-family $Y^V$ as follows.
Let $Y^V$ denote the space of all $\V$-maps 
$V\to Y$ with the compact-open topology. 
Then we obtain the 
projection 
\begin{equation}
Y^V \to \bar Y
\end{equation}
which sends a $\V$-map $a\from V \to Y$ 
to the point $p_Ya(v_0)\in \bar Y$, where $v_0\in V$ is 
an arbitrary point of 
the fibre $V$. The pre-image of a point $b\in \bar Y$
under this projection is equal to $\ff \hom(V,Y_b) = Y_b^V$.
Thus $Y^V \to \bar Y$ 
is a $(\V,\ff^V)$-family.

We recall some properties of the function spaces
$X^V$ and $W^V=\hom_\V(V,W)$.
By assumption $(*)$ the images under the fibre functor $\ff$ of 
all the objects in $\V$ are locally compact 
second-countable and Hausdorff. This implies that they 
are metrizable. 
Since every object $V$ is metrizable, 
the automorphism group $\Aut(V)$
is a topological group with the compact-open topology
induced by the inclusion $\Aut(V) \subset \hom_\V(V,V)$.
% ref:
% Arens, Amer J Math 68 (1946) 593--610.

Given a $\V$-map $f\from X \to Y$, let $f^V$ denote the map
$f^V\from X^V \to Y^V$ defined by composing with $f$ the
$\V$-maps $V\to X$. 
Moreover for $\varphi\from V \to W$ in $\V$
let $X^\varphi\from X^W \to X^V$
be the map induced by $\varphi$.
Let us recall that given a $\V$-map
$f\from Z\times V \to X$ the adjoint $\hat f$ of $f$ is the function
$\hat f\from Z\to X^V$ defined by $\hat f(z)(v) = f(z,v)$ for every $z\in Z$
and $v\in V$.

\begin{lemma} \label{lemma:adjoint}
For every space $Z$
a $\V$-function $f\from Z\times V \to X$ is continuous if and 
only if the adjoint 
$\hat f: Z \to X^V$ is continuous.
\end{lemma}
\begin{proof}
Because $V$ is locally compact and Hausdorff, a function
$f\from Z \times V \to X$ is continuous if and only if 
its adjoint $\hat f: Z \to \mathrm{Map}(V,X)$ is continuous,
where $\mathrm{Map}(V,X)$ denotes the space of all (not necessarily
controlled by $\V$) maps from $V$ to $X$ with the compact-open topology.
But a function $f\from Z\times V \to X$ is a $\V$-map if and only 
if its adjoint sends $Z$ into the subspace 
$X^V \subset \mathrm{Map}(V,X)$. This implies the lemma.
\end{proof}

\begin{lemma} \label{coro:adjoint}
The evaluation map
$X^V\times V \to X$ which sends $(g,v)$ to $g(v)$ is continuous.
Moreover, if $f\from X\to Y$
is continuous then the induced function $f^V\from X^V \to Y^V$
is continuous. Also, given $\varphi\from V \to W$ 
the induced map $X^\varphi\from X^W \to X^V$ is continuous.
\end{lemma}
\begin{proof}
The evaluation map is the adjoint of the identity of $X^V$,
hence continuous by lemma \ref{lemma:adjoint}.
The evaluation at $v\in V$ is the restriction of 
the evaluation to the subspace $X^V\times\{v\}$ of $X^V\times V$
and hence continuous. The projection $p\from X^V \to \bar X$
is the composition of the evaluation at any $v_0\in V$ with 
the projection $p_X$ and hence continuous.
The induced function $f^V$ is the adjoint of the composition
$X^V\times V \to X \to Y$ where the first arrow is the evaluation
and the second is $f$.
To see that $X^\varphi$ is continuous it is enough to see
that it is the adjoint of the composition $X^W\times V \to X^W \times W \to X$,
where the first arrow is $1\times \varphi$ while the second
is the evaluation.
\end{proof}

The next result is a crucial observation which shows that the function
space $X^V$ of a stratified bundle $X$ is again a stratified bundle.

\begin{theo}
\label{theo:principal}
Let $\V$ be a NKC category with fibre functor $\ff$.
If $X$ is a 
% $(\V,\ff)$-complex,
$(\V,\ff)$-stratified bundle
and $V$ an object in $\V$
then the function space 
$X^V$ is also a %$(\V,\ff^V)$-complex,
$(\V,\ff^V)$-stratified bundle.
\end{theo}
\begin{proof}
If $X$ is a $(\V,\ff)$-complex then 
this follows from corollary \ref{coro:main} below.
In the general case one can apply a finite number of times
corollary 
\ref{coro:nkc} and lemma \ref{lemma:nnep2}.
\end{proof}

Let $\top_{\bar X}$ be the category of spaces over $\bar X$;
objects are maps $X\to \bar X$.
Given a $\V$-family $X$, 
the function spaces $X^V$ for every object $V$ in $\V$
yield a functor 
$X^\circ$
which carries $V$ to $X^\circ(V) = X^V$.
Moreover, given a $\V$-map $f\from X \to Y$ there is a natural transformation
$f^\circ: X^\circ \to Y^\circ$ defined by $f^V\from  X^V \to Y^V$.
Let $\Vdiag_{\bar X}$ denote the diagram category 
in which 
the objects
are functors $\V^\op \to \top_{\bar X}$ and the morphisms are the 
natural transformations between functors.
The operator $(-)^\circ$ 
sends a $\V$-family $X$ to the $\V^\op$-diagram
$X^\circ$
and a  $\V$-map  $f$ 
to  the natural transformation $f^\circ$.
Thus ${(-)}^\circ$ is a functor
\begin{equation*}
{(-)}^\circ\from \Vtop_{\bar X} \to \Vdiag_{\bar X}.
\end{equation*}

The $\V$-diagram $X^\circ$ is the generalization
of the concept of principal bundle.
In fact, 
if $X\to \bar X$ is a fibre bundle with structure group $G$ and 
% locally compact second-countable Hausdorff
fibre $V$
and $\V=G \to \top$ is given by the $G$-space $V$
then  
the diagram $X^\circ\from \V^\op \to \top$
coincides with the $G$-space $X_G$
given by the principal bundle $X_G\to\bar X$ associated
to $X$, see \ref{eq:steenrod}.

For  a topological enriched category $\V$ there is 
a naturally associated  functor
$\hm\from \V \to \V^\op$-$\Diag$ which sends an object $W$ of 
$\V$  to the $\V^\op$-diagram $\hom_\V(-,W) = W^\circ$.
In case $X$ is a $\V$-stratified bundle 
and $(\V,\ff)$ is a NKC structure category
then by theorem \ref{theo:principal} 
the function space $X^V$ is not only a space over $\bar X$ but
a $(\V,\ff^V)$-stratified bundle.
Moreover this holds for every object $V\in \V$
in a compatible way.
Therefore the diagram $X^\circ$ given by stratified bundles $X^V$
leads to the concept of a 
\emph{stratified bundle diagram}.
For this consider the definition of a stratified bundle
in \ref{defi:stratifiedvectorbundle}
and replace the category $\top$ by the category of diagrams
in $\top$. This yields a notion of a stratified bundle
diagram. 
If the stratification is given by the skeletal filtration 
of the base space $\bar X$
then a $(\V,\hm)$-stratified bundle diagram is the same 
as a ``free CW-complex'' in the sense of Davis--L\"uck
\cite{davislueck}.

\begin{defi}\label{defi:principal}
If $\V$ is a topologically enriched category,
then  
a $(\V,\hm)$-stratified bundle diagram 
in the category $\V^\op$-$\Diag_{\bar X}$
is termed \emph{principal stratified bundle diagram}.
This is a $\V^\op$-diagram of stratified $\V$-bundles 
over a stratified space $\bar X$.
For example $X^\circ\from \V^\op\to \top_{\bar X}$ is such a principal
stratified bundle diagram.
\end{defi}

If $\V$ has a single object $V$ and $\hom_\V(V,V)$ is a topological
group,
then a $(\V,\hm)$-stratified bundle  diagram
is nothing but
a principal $G$-bundle with structure group $G=\hom_\V(V,V)$ over $\bar X$,
see theorem \ref{theo:bundle}.

Now consider a principal stratified bundle  diagram
$P\from \V^{\op} \to \top_{\bar X}$ over $\bar{X}$ 
and a fibre functor $\ff\from \V \to \top$.
We can build a $(\V,\ff)$-stratified bundle $X\to \bar X$ in $\top$
as the coend of the two functors $P$ and $\ff$:
\begin{equation*}
X = P\otimes_\V \ff = \int^{V\in \V} P(V)\times \ff(V).
\end{equation*}
We recall the coend construction: the space $P\otimes_{\V} \ff$ is 
the quotient space of the coproduct
\begin{equation*}
\coprod_{V\in \V} P(V)\times \ff(V)
\end{equation*}
with the identification 
$(P(\alpha)(x),y) \sim (x,\ff(\alpha)(y) )$
for every morphism $\alpha\in \hom_\V(V,W)$,
$x\in P(W)$ and $y\in \ff(V)$. There is a well-defined 
map $X \to \bar X$ and it is not difficult to show:

\begin{lemma}\label{lemma:coend}
$P\otimes_\V \ff$ is a $(\V,\ff)$-family.
\end{lemma}
\begin{proof}
See e.g. \cite{huse}, 
page 44, for the case of a groupoid $\V$. 
It suffices to show it in the case $P$ is a $(\V,\hm)$-point, that is 
a diagram $P$ over a point. Since $P$ is principal,
there is an object $W$ of $\V$ such that $P(V)=\hom_{\V}(V,W)$
for every $V\in \V$, or, equivalently, $P = W^\circ$. 
If $g\from V \to V'$ is a morphism in $\V$, 
then $P(g)\from W^{V'} \to W^V$ is defined by
$P(g)(\alpha) = \alpha g$. Thus the coend $W^\circ \otimes_\V \ff$
is defined as 
\begin{equation}
\left(\coprod_{V\in \V} W^V \times \ff(V) \right)/\sim, 
\end{equation}
where $(\alpha g, x) \sim (g, \ff(\alpha)(x) )$
for every $g\in \hom_\V(V,V')$, $\alpha\in W^V$ and 
$x\in \ff(V)$.
Now it is easy to see that for every $W$
\[
W^\circ \otimes_\V \ff = \ff(W)
\]
and that given a morphism $\alpha^\circ \from W^\circ\otimes_{\V} \ff \to 
{W'}^\circ \otimes \ff$ the following diagram commutes:
\ldiag{diag:nota}{%
W^\circ \otimes_\V \ff \ar[r]^{\alpha^\circ} \ar[d]^{=} & 
{W'}^\circ \otimes_\V \ff \ar[d]^{=} \\
\ff(W) \ar[r]^{\ff(\alpha)} & \ff(W') \\
}%
\end{proof}

\begin{theo}\label{theo:principalbundle} %\label{theo:last}
Let $\V$ be a structure category with a faithful fibre functor
$\ff$ so that $\V$ is topological enriched by the compact
open topology.
Let $(\V,\ff)$ be a NKC structure category and $X$ a 
$(\V,\ff)$-stratified  bundle.  
Then the associated diagram $X^\circ$
is a principal stratified bundle diagram
with fibre functor $\hm$.
Conversely, if $P$ is a principal
$(\V,\hm)$-stratified bundle diagram and $\ff$ a fibre functor, 
then the coend $P\otimes_\V \ff$ is a stratified 
bundle in $\top$ with structure category $\V$ and fibre functor $\ff$.
Moreover 
\[
X^\circ \otimes_\V \ff = X.
\]
\end{theo}
This result generalizes the classical equation 
\[
X_G\times_G V = X
\]
for the principal bundle $X_G$ associated to $X$ in \ref{diag:5:1}.
This construction can be generalized. Let $\varphi\from G \to H$
be a continuous homomorphism between topological groups
and let $G$ ($H$) be a topological transformation group for 
$V$ ($W$). Then a \emph{$(G,V)$-bundle} $X$ (i.e. a bundle with fibre 
$V$ and structure group $G$) yields
the principal bundle $X_G$ which in turn yields
the \emph{$\varphi$-associated bundle}
\[
\varphi_\#(X) = X_G \times_G \varphi^*W \to \bar X.
\]
Here $\varphi^*W$ is the $G$-space with the $G$-action induced by
$\varphi$, that is $g\cdot w = \varphi(g)\cdot w$ for $g\in G$
and $w\in W$.
The associated bundle $\varphi_\#(X)$ has fibre $W$ and structure
group $H$. Hence $\varphi$ 
induces the functor
\newcommand{\bundles}[1]{\mbox{$(#1)$-${\mathbf{Bundles}}_{\bar X}$}}
\begin{equation}
\label{eq:varphi}
\varphi_\#\from \bundles{G,V} \to \bundles{H,W}.
\end{equation}
We now generalize this functor for stratified bundles.

Let $(\V,\ff)$ and $(\W,\gg)$ be structure categories
with faithful fibre functor $\ff$, $\gg$ so that 
$\V$ and $\W$ have the compact open topology.
Let 
\[
\varphi\from \V \to \W
\]
be a continuous functor.
If $(\V,\ff)$ is a NKC structure category and $X$ a $(\V,\ff)$-stratified
bundle then the principal stratified bundle
diagram $X^\circ$ is defined and we obtain the 
\emph{$\varphi$-associated stratified bundle}
\[
\varphi_\#(X) = X^\circ \otimes_\V (G\varphi) \to \bar X.
\]
Here $\varphi_\#(X)$ is a $(\W,\gg)$-stratified bundle.
Hence $\varphi$ induces the functor
\[
\varphi_\#\from (\V,\ff)\mbox{-}\Stra_{\bar X} \to
(\W,\gg)\mbox{-}\Stra_{\bar X}.
\]
\begin{proof}[Proof of theorem \ref{theo:principalbundle}]
We have seen that, as a consequence of theorem \ref{theo:principal},
for each $V\in \V$ the function space $X^V$ is obtained by
gluing $W^V$-bundles $M_i^V$ along with maps  
$h_i^V\from A_i^V\subset M_i \to X_{i-1}^V$,
for suitable objects $W\in \V$.
In other words, the diagram $X^\circ$ is obtained by gluing
the $\V$-diagrams $M_i^\circ$ (which are  $(\V,\hm)$-bundles)
along with $\V$-diagram maps $h_i^\circ\from A_i^V \to X_{i-1}^\circ$.
Thus $X^\circ$ is a $(\V,\hm)$-stratified bundle, that is 
a $\V$-stratified bundle with diagram-fibre functor $\hm$.
By definition this is a principal stratified bundle diagram.

Conversely, 
let $P$ be a principal stratified bundle diagram 
and $\ff$ a fibre functor.
It is not assumed that $\ff$ is faithful, but that the spaces $\ff V$
are second-countable locally compact Hausdorff spaces,
see assumption $(*)$.
Hence, if the $i$-skeleton of $P$ is obtained by
attaching a $(\V,\hm)$-bundle  $\hat M_i$ to $P_{i-1}$ via 
a $\V$-diagram map $\hat h_i\from \hat A_i\subset \hat M_i \to P_{i-1}$, i.e.
$P_i = \hat M_{i} \cup_{\hat h_i} P_{i-1}$,
the coends are corners of a push-out diagram in $\Vtop$:
\ndiag{%
\hat A_i\otimes_\V \ff \ar[r]^{\hat h_i\times_\V \ff} \ar[d] 
& 
P_{i-1}\otimes_\V \ff \ar[d] \\
\hat M_i\otimes_{\V} \ff \ar[r] 
& P_i\otimes_{\V}\ff\\
}%
Thus $X$ has a filtration given by
$X_i = P_i\otimes_\V \ff$, for all $i\geq 0$,
and each $X_i$ is obtained by attaching $M_i = \hat M_i \otimes_\V \ff$
to $X_{i-1}$ via the map 
\[
h_i= \hat h_i \otimes_\V \ff \from A_i = \hat A_i\otimes_\V \ff \to
X_{i-1}.
\]
Moreover, since 
\[
\lim_{i\geq 0}  \left( P_i\otimes_\V  \ff \right) = 
\left( \lim_{i\geq 0}   P_i \right) \otimes_\V  \ff,  
\]
$X$ is equal to the colimit $X=\lim_{i\geq 0} X_i$.
To show that it is a $(\V,\ff)$-stratified  bundle we need only to 
show that the attaching maps are $\V$-maps. But $\hat h_i$ are 
$(\V,\hm)$-maps, and by diagram \ref{diag:nota} the coends
of such maps are $(\V,\ff)$-maps.
The proof is hence complete.
\end{proof}

\section{Pull-back of a push-out $\V$-family}
\label{section:pullback}
In this section we complete the proof of 
the pull-back theorem \ref{theo:pullback}.

Consider a $\V$-family $Y'$, a $\V$-CW-pair $(M',A')$
and a $\V$-map $h'\from A' \to Y'$. Then we can 
glue $M'$ to $Y'$ and obtain a $\V$-family $X'$, as in the following diagram.
\ldiag{diag:fri1}{%
A' \ar[rrr]^{h'} 
\ar[ddd]^{j_{A'}}
\ar[rd]^{p_{A'}}
&
&
&
Y' \ar[ddd]^{j_{Y'}}
\ar[dl]^{p_{Y'}}
\\
&
\bar A'
\ar[r]^{\bar h'}
\ar[d]^{j_{\bar A'}}
&
\bar Y' \ar[d]^{j_{\bar Y'}} 
\\
&
\bar M' 
\ar[r]^{\bar \hh'}
&
\bar X' \\
M' \ar[rrr]^{\hh'} \ar[ur]^{p_{M'}}
& &  &
X' \ar[ul]^{p_{X'}}\\
}%

Now consider a space $\bar Y$, 
a CW-pair $(\bar M, \bar A)$
and a map $\bar h\from \bar A \to \bar Y$.
Let $\bar X$ be the push-out space 
$\bar X = \bar Y \cup_{\bar h} \bar M$.
Assume that maps $\tilde f_A$, $\tilde f$,
$\bar f_{\bar Y}$ and $\bar f$ are given, such 
that the following diagram commutes.
\ldiag{diag:fri2}{%
\bar A' \ar[rrr]^{\bar h'}
\ar[ddd]^{j_{A'}}
&
&
&
\bar Y'
\ar[ddd]^{j_{\bar Y'}}
\\
&
\bar A
\ar[ul]^{\tilde f_A}
\ar[r]^{\bar h}
\ar[d]^{j_{\bar A}}
&
\bar Y
\ar[d]^{j_{\bar Y}}
\ar[ur]^{\bar f_{\bar Y}}
\\
& 
\bar M
\ar[r]^{\bar \hh}
\ar[dl]^{\tilde f}
&
\bar X 
\ar[dr]^{\bar f}
\\
\bar M'
\ar[rrr]_{\bar \hh'}
&&&
\bar X'
\\
}%

Then the pull-back families  $X=\bar f^* X'$, $Y=\bar f_{\bar Y}^* Y'$,
$A= \tilde f_{A}^* A'$ and 
$M = \tilde f^*M'$ are defined.
Moreover, we can define the maps $h=(\bar h,h')\from A \to Y$,
$j_Y=(j_{\bar Y},j_{Y'})\from Y \to X$,
$j_A=(j_{\bar A},j_{A'})\from A \to M$
and $\hh=(\bar \hh, \hh')\from M \to X$.
They fit into the following diagram.
\ldiag{diag:fri3}{%
\bar A \times A' 
&
\ar@{ >->}[l]
A \ar[r]^{h}
\ar[d]^{j_A}
&
Y 
\ar@{ >->}[r]
\ar[d]^{j_Y}
&
\bar Y\times Y' 
\\
\bar M \times M' 
&
\ar@{ >->}[l]
M
\ar[r]^{\hh}
&
X 
\ar@{ >->}[r] &
\bar X \times X'.
}%

The aim of this section is to prove the following proposition.
\begin{propo}\label{propo:theo:pushout2}
Assume that $\bar M$ and $\bar M'$ are locally finite
and countable CW-complexes 
and that  
$\bar f(\bar X \minus \bar Y) \subset \bar X' \minus \bar Y'$.
Then the 
pull-back family $X=\bar f^*X'$
is obtained by attaching the pull-back $\V$-family
$M=\tilde f^* M'$ to 
$Y=\bar f_{\bar Y}^*Y'$ via the $\V$-map $h\from A\subset M \to Y$
induced by $(\bar h, h')$, i.e. the middle square in diagram
\ref{diag:fri3} is a push-out.
\end{propo}

The proof of proposition \ref{propo:theo:pushout2} will be 
the content of the rest of the section.

\begin{lemma}\label{lemma:step1}
The maps in \ref{diag:fri3} are well-defined and the diagram commutes.
\end{lemma}
\begin{proof}
Consider $(\bar a, a')\in A$. Then $\tilde f_A(\bar a) = p_{A'}(a')$,
hence 
\begin{equation*}
\bar f_{\bar Y} \bar h (\bar a) =
\bar h' \tilde f_A (\bar a) = \bar h'
 p_{A'}(a') 
 =
 p_{Y'} h'(a').
\end{equation*}
Hence 
$(\bar h(\bar a), h'(a') ) \in Y$.
On the other hand, 
if $(\bar y, y' ) \in Y$, then
$\bar f_{\bar Y}(\bar y) = p_{Y'}(y')$, hence
\begin{equation*}
\bar f j_{\bar Y} (\bar y) =
j_{\bar Y'} \bar f_{\bar Y} (\bar y) =
j_{\bar Y'} p_{Y'}(y') =
p_{X'} j_{Y'} (y').
\end{equation*}
That is,
$j_Y( \bar y, y') = (j_{\bar Y} (\bar y), j_{Y'}(y') ) \in X$. 
The same argument can be applied literally to $j_A$ and 
$\hh$.
The diagram commutes since it is the restriction of a commutative
diagram on the Cartesian products.
\end{proof}

\begin{lemma}\label{lemma:step2}
The maps $j_A\from A \to M$  and $j_Y\from Y \to X$ 
are  a closed inclusions of $\V$-families.
\end{lemma}
\begin{proof}
Consider the following diagram. The left square is a pull-back
(by definition) and the right square is a pull-back
(the pair $(M',A')$ is a $\V$-CW-pair).
\ldiag{diag:fri4}{%
A \ar[r]
\ar[d]
&
A'
\ar[r]\ar[d]
&
M' \ar[d]\\
\bar A \ar[r]
&
\bar A' \ar[r]
&
\bar M' \\
}%
Now consider this diagram.
\ndiag{%
A \ar[r] \ar[d] 
&
M \ar[r] \ar[d]
&
M' \ar[d]
\\
\bar A \ar[r] 
&
\bar M \ar[r]
& \bar M'\\
}%
The composition of the two squares is equal to the composition
of the two squares in diagram \ref{diag:fri4},
hence it is a pull-back. The right square is a pull-back by definition,
hence the left square is a pull-back.
(see e.g. exercise 8 of \cite{macl}, page 72).
Since $\bar A \subset \bar M$ is a closed inclusion,
$j_A$ is a closed inclusion. The proof is the same for $j_Y$.
\end{proof}

Now consider the maps $q=\hh \sqcup j_{Y}$, 
$\bar q = \bar \hh \sqcup j_{\bar Y}$,
$\bar q' = \bar \hh' \sqcup j_{\bar Y'}$
and 
$q' = \hh' \sqcup j_{Y'}$.
They can be arranged in the following diagram.
\ldiag{diag:fri5}{%
M\sqcup Y \ar[rrr] \ar[ddd]^q \ar[rd]
&&&
M' \sqcup Y' \ar[ddd]^{q'} 
\ar[dl]
\\
&
\bar M \sqcup \bar Y 
\ar[r]\ar[d]^{\bar q}
&
\bar M' \sqcup \bar Y' \ar[d]^{\bar q'}
\\
&
\bar X \ar[r] & \bar X' \\
X \ar[rrr] \ar[ur] &&& X' \ar[ul] \\
}%

By definition $\bar q$, $\bar q'$ and $\bar q'$ are 
quotient maps. We want to show that under suitable conditions 
$q$ is a quotient map.
If this is the case, then proposition \ref{propo:theo:pushout2} is proved,
since $q$ is exactly the projection defining the topology of 
the push-out.

\begin{lemma}\label{lemma:step4}
Assume that $\bar f(\bar X \minus \bar Y) \subset \bar X' \minus \bar Y'$.
Then $q$ is onto, $q|Y$  and $q|(M\minus A)$ are mono.
\end{lemma}
\begin{proof}
Let $x=(\bar x, x') \in X$, so that 
then $\bar f ( \bar x ) =  p_{X'} (x')$.
If $\bar x \in \bar Y$, then $\bar f (\bar x) \in \bar Y'$,
therefore $x' \in p_{X'}^{-1} \bar Y' = Y'$.
Thus there is $y = (\bar y, y') \in Y $ such that 
$q(y) = x$.
On the other hand, 
if $\bar x\not\in \bar Y$, since assumption 
$\bar f \bar x \not\in \bar Y'$,
necessarily $x' \not\in  Y'$
hence there is a unique $m'\in M'\minus A'$
such that $q'(m') = x'$.
For the same reason there is a unique $\bar m\in \bar M\minus \bar A$
such that $\bar q(\bar m) = \bar x$.
Now consider the since $\bar q'$ is mono in $\bar M' \minus A'$,
the chain of equalities
\begin{equation*}
\bar q' \tilde f (\bar m) 
=
\bar f \bar q(\bar m)  = \bar f (\bar x) 
= p_{X'}( x') 
=
p_{X'} q' (m')
=
\bar q' p_{M'} (m')
\end{equation*}
implies that $\tilde f (\bar m) = p_{M'}(m')$,
hence that 
$(\bar m, m') \in M$.
We have shown that $q$ is onto.
The restriction of $q$ to $Y$ 
is mono 
since 
$j_{\bar Y} 
\times j_{Y'}$ is mono.
To see that the restriction of $q$ to 
$M\minus A$ is mono, consider 
that if $(\bar m, m')\in M$ then 
$\bar m \not\in \bar A$ implies 
$m'\not \in A'$,
hence 
$M \minus  A \subset (\bar M\minus \bar A)\times (M'\minus A')$.
Now, $\bar q$ and $q'$ are mono when restricted 
to $\bar M \minus \bar A$ and 
$M'\minus A'$ respectively, hence $q|(M\minus A)$ is mono.
\end{proof}

From now on, we will assume that the condition of lemma \ref{lemma:step4}
is fulfilled.

\begin{lemma}
\label{lemma:step7}
The restriction $q|(M\minus A)\from M\minus A \to X\minus Y$
is a homeomorphism.
\end{lemma}
\begin{proof}
The map $q|(M\minus A)$ is bijective and continuous. We need to show that 
it is open. Consider an open set $U$ in $M\minus A$,
and a point $x=(\bar x, x') \in U$.  
There are neighborhood $U_{\bar M} \subset \bar M \minus \bar A$
and $U_{M'} \subset M'\minus A'$ 
such that 
\[
x \in ( U_{\bar M} \times U_{M'} ) \cap M \subset U.
\]
Since $\bar q$ and $q'$ are homeomorphisms of $\bar M \minus \bar A$
and $M'\minus A'$ onto their images,
the sets $U_{\bar X} = \bar q U_{\bar M}$
and 
$U_{X'} = q' U_{M'}$ 
are open subsets of $\bar X \minus \bar Y$
and $X' \minus Y'$.
Hence $U_{\bar X}\times U_{X'} \cap X$ 
is open in $X$ and contained in $X\minus Y$. Moreover, since $q$ 
is mono in $M\minus A$, we have that 
$q \left[ ( U_{\bar M} \times U_{M'} ) \cap M \right]
= U_{\bar X} \times U_{X'} \cap X $,
hence $U_{\bar X} \times U_{X'} \cap X \subset q(U)$.
But $x$ is arbitrary, therefore $q(U)$ is open
and so the map $q|M\minus A$ is open.
\end{proof}

\begin{lemma}
\label{lemma:step6}
A subset $S\subset M\sqcup Y$ is saturated (i.e. $q^{-1}q(S) = S$ )
if and only if
$S\cap A = h^{-1}(S\cap Y)$.
\end{lemma}

\begin{lemma}\label{lemma:longnonce}
Assume that 
$\bar M$ and $M'$ are metrizable.
%(if and only if it is regular and has a basis that can be decomposed into 
% an at most countable collection of nbd-finite families) . 
% nagata page 194 dugundji (check fritsch piccinini)
% $\bar A$ and $A'$ are locally compact 
% or that they have a countable number of cells.
If $U_{\bar M} \subset \bar M$ and 
$U_{A'} \subset A'$ are open sets with compact closure such that 
\[
(\overline{U_{\bar M}} \times  \overline{U_{A'}}) \cap M \subset  q^{-1}U, 
\]
then there is an open subset 
$U_{M'} \subset M'$ with compact closure such that 
\begin{gather*}
U_{M'} \cap A' = U_{A'}\\
(\overline{U_{\bar M}} \times  \overline{U_{M'}}) \cap M \subset  q^{-1}U.
\end{gather*}
In the same way,
if $U_{\bar A} \subset \bar A$ and $U_{M'}\subset M'$ 
are open sets with compact closure such that 
\[
(\overline{U_{\bar A}} \times  \overline{U_{M'}}) \cap M \subset  q^{-1}U, 
\]
then there is an open subset 
$U_{\bar M} \subset \bar M$ with compact closure such that 
\begin{gather*}
U_{\bar M} \cap \bar A = U_{\bar A}\\
(\overline{U_{\bar M}} \times  \overline{U_{M'}}) \cap M \subset  q^{-1}U.
\end{gather*}
\end{lemma}
\begin{proof}
By assumption the set $q^{-1}U\cap M$ is open in $M\subset \bar M \times M'$, 
therefore there exists an open set 
\newcommand{\tU}{\widetilde{q^{-1}U}} 
$\widetilde{q^{-1}U} \subset \bar M \times M'$ such that 
$\tU \cap M = q^{-1}U \cap M$.
We need to find the open set $U_{M'}\subset M'$ with the desired properties.
Consider the following function 
$\eta\from \bar M \times M' \to \R$
defined by
\begin{equation}\label{eq:function}
\eta(x) = d(x,{\tU}^c) + d(x,M),
\end{equation}
where $\tU^c$ is the complement of $\tU$ in $\bar M \times M'$ 
and $d$ is the distance given by a metric on $\bar M \times M'$.
Consider  a point $x$ in the compact 
closure $\overline{U_{\bar M}} \times \overline{U_{A'}}$.
If $x\in M$ as well, then by the assumption
$x\in \tU$, and since $\tU$ is open this implies $\eta(x)\geq d(x,{\tU}^c)>0$.
On the other hand, since $M$ is closed in $\bar M \times M'$,
if $x\not\in M$ then $\eta(x)\geq d(x,M) >0$.
Thus $\eta(x)>0$ for every 
$x\in (\overline{U_{\bar M}} \times \overline{U_{A'}} )$,
and so there is an open set $U_{M'} \subset M'$ with compact closure
such that $U_{M'} \cap A' = U_{A'}$
and 
\[
x\in \overline{U_{\bar M}} \times \overline{U_{M'}} \implies \eta(x)>0.
\]
This means that if 
$x\in (\overline{U_{\bar M}} \times \overline{U_{M'}}) \cap M$ then 
$d(x,\tU^c)>0$, thus $x\in \tU$, and this implies that 
\[
(\overline{U_{\bar M}} \times \overline{U_{M'}}) \cap M \subset q^{-1}U.
\]
For the second part of the lemma, 
the proof is exactly the same,
it is only necessary to exchange the roles of $\bar A$ and $M'$.
\end{proof}

\begin{lemma}\label{lemma:a1}
Let $X$ be a $\V$-bundle over 
a locally finite and countable CW-complex $\bar X$.
Then $X$ is metrizable. Moreover, 
every open subset $O\subset X$ 
is the union of an ascending sequence of open subsets $O_n \subset X$
with compact closure $\overline{O_n} \subset O_{n+1}$.
\end{lemma}
\begin{proof}
Since $\bar X$ is locally finite and countable,
it is metrizable and $2^\circ$ countable. 
Since the fibres are $2^\circ$ countable by assumption,
$X$ is 
$2^\circ$ countable,
since it is obtained by attaching a countable number 
of $\V$-cells $V\times D^n$, which are $2^\circ$ countable.
Now, 
$\bar X$  
and the $\V$-cells 
$V\times D^n$
are completely regular,
hence $X$ is regular;
therefore by
Urysohn metrization theorem $X$ is metrizable. 
Let $d$ denote its metric. 
Now consider an ascending chain of finite subcomplexes 
$\bar X_n \subset \bar X$  such that $\bar X_n$ 
is contained in the interior of $\bar X_{n+1}$ 
and $\bar X=\bigcup_n \bar X_n$ (such a sequence exists 
since $\bar X$ is locally finite and countable).
Using the fact that the fibres are locally compact and 
second countable, it is possible to show that 
$X$ can be written as the union of an ascending sequence of 
subspaces  $X_n \subset X_{n+1}\subset \dots$
where 
each $X_n$ is compact and is contained in the interior 
$\mathring{X}_{n+1}$.
For every $n$ define the following open set
\[
O_n = \{ 
x\in X \st
d(x,O^c)>\frac{1}{n} 
\} \bigcap \mathring{X_n},
\]
where $O^c$ denotes the complement of $O$.
If $x$ is an element of the closure $\overline {O_n}$ (closure of $O_n$ in $X$),
then $x \in O_{n+1} \bigcap X_n \subset X_n$,
hence $\overline{O_n}$ is a closed subset of the compact space $X_n$, hence
$\overline{O_n}$ is compact. Moreover, 
since $X_n \subset \mathring{X}_{n+1}$, we have that 
$\overline{O_n} \subset O_{n+1}$ as claimed.
We need to show that $O = \bigcup_n O_n$:
if $x\in O$, then $d(x,O^c)>0$, therefore 
there is $n_1\geq 1$ such that $d(x,O^c)>\frac{1}{n}$
for every $n\geq n_1$. 
Since $\bigcup_n X_n = X$, there is $n_2\geq n_1$ such that 
$x\in X_{n_2}\subset \mathring{X}_{n_2+1}$.
But these conditions imply that 
$x\in O_{n}$ for $n=n_2+1$. Thus $O=\bigcup  O_n$.
\end{proof}

\begin{lemma}\label{lemma:nonce}
Assume that 
$\bar M$ and $M'$ are metrizable.
If $U_{\bar A} \subset \bar A$ and 
$U_{A'} \subset A'$ are open sets such that 
\[
({U_{\bar A}} \times  {U_{A'}}) \cap M \subset  q^{-1}U, 
\]
then there are open sets $U_{\bar M} \subset \bar M$ and 
$U_{M'} \subset M'$ such that 
\begin{gather*}
U_{\bar M} \cap \bar A = U_{\bar A}\\
U_{M'} \cap A' = U_{A'}\\
(U_{\bar M} \times  U_{M'}) \cap M \subset  q^{-1}U.
\end{gather*}
\end{lemma}
\begin{proof}
By lemma \ref{lemma:a1},
it is possible to find increasing sequences 
$U_{\bar A}^k \subset U_{\bar A}$ and
$U_{A'}^k \subset U_{A'}$ with compact closures
and such that 
\begin{gather*}
\bigcup_{k\geq 1} U_{\bar A}^k = U_{\bar A} \\
\bigcup_{k\geq 1} U_{A'}^k = U_{A'}.
\end{gather*}
Now, by applying lemma \ref{lemma:longnonce} twice  
it is possible to define compact subsets $U_{\bar M}^1$ 
and $U_{M'}^1$
such that 
\begin{gather*}
(\overline{U_{\bar M}^1} \times  \overline{U_{M'}^1} ) \cap M \subset q^{-1}U.
\end{gather*}
By induction, we will show that it is possible to find 
two sequences of 
open sets $U^k_{\bar M}$ 
and $U_{M'}^{k}$ with compact closures, with the property that 
\begin{gather*}
U^k_{\bar M} \cap \bar A = U_{\bar A}^k \\
U^{k}_{M'} \cap A' = U_{A'}^{k} 
\end{gather*}
for every $k$, $k' \geq 1$ and 
\begin{equation}\label{eq:nonce}
( \overline{ U_{\bar M}^k} \times \overline{ U_{M'}^{k'}} ) \cap M \subset q^{-1}U
\end{equation}
for every $k$, $k' \geq 1$.
Assume that the sequences are defined for $j\leq k-1$.
By applying $k-1$ times lemma \ref{lemma:longnonce} 
and taking the intersection of the resulting open sets,
we can show that there is $U_{\bar M}^k$ 
such that 
$U^k_{\bar M} \cap \bar A = U_{\bar A}^k$
and 
\[
( \overline{ U_{\bar M}^k} \times \overline{ U_{M'}^{k'}} ) \cap M \subset q^{-1}U
\]
for every $k'\leq k-1$.
Now we can apply $k$ times lemma \ref{lemma:longnonce} 
and take the intersection of the resulting open sets,
to finally find and open set with compact closure $U_{M'}^k$ such that 
$U^k_{M'} \cap A' = U_{A'}^k$ and 
\[
( \overline{ U_{\bar M}^j} \times \overline{ U_{M'}^{k}} ) \cap M \subset q^{-1}U
\]
for every $j\leq k$. Thus equation \ref{eq:nonce}
holds for every $k$, $k'\geq 1$.
But this implies that the open sets 
\begin{gather*}
U_{\bar M} = \bigcup_{k\geq 1} U^k_{\bar M} \\
U_{M'} = \bigcup_{k\geq 1} U^k_{M'}
\end{gather*}
have the desired property.
\end{proof}

\begin{lemma}\label{lemma:step8}
If $\bar M$ and $\bar M' $ are locally finite and countable, then
the map $q$ is a quotient map.
\end{lemma}
\begin{proof}
Consider a subset $U\subset X$ 
such that $q^{-1}U$ is open in $M\sqcup Y$. 
We want to show that $U$ is open.
Let $x\in U \minus Y$. 
Since $q^{-1}U\cap (M\minus A)$ is open,
by lemma \ref{lemma:step7}
there is an open neighborhood of $x$ in $X\minus Y$ 
contained in $U$.
If $U\cap Y=\emptyset$ then we have proved that $U$ is open.
Otherwise, let $x\in Y \cap U$. 
Since $q^{-1}U$ is open, $q^{-1}U \cap Y$ is open
and contains $y\in q^{-1}(x)\cap Y$.
Therefore there are two open sets 
$U_{\bar Y}$ and $U_{Y'}$ 
such that 
\begin{equation*}
y \in (U_{\bar Y} \times U_{Y'} )\cap Y \subset q^{-1}Y \cap Y.
\end{equation*}
Let 
\begin{gather*}
U_{\bar A} = \bar h^{-1} U_{\bar Y}\\
U_{A'} = {h'}^{-1} U_{Y'}
\end{gather*}
Since $\bar h$ and $h'$ are continuous,
they are open subsets of $\bar A$ and $A'$ respectively.
By lemma \ref{lemma:nonce}
there are open  sets
$U_{\bar M}  \subset \bar M$ 
and $U_{M'} \subset M'$ 
such that 
\begin{gather*}
U_{\bar M} \cap \bar A = U_{\bar A}= \bar h^{-1} U_{\bar Y}\\
U_{M'} \cap A' = U_{A'}= {h'}^{-1} U_{Y'},
\end{gather*}
and 
\[
(U_{\bar M} \times U_{M'})\cap M \subset q^{-1}U.
\]
Consider the sets
\begin{gather*}
U_{\bar X} = \bar q ( U_{\bar M} ) \cup U_{\bar Y} \subset \bar X \\
U_{X'} = q' ( U_{M'} ) \cup U_{Y'} \subset X'.
\end{gather*}
Since $\bar q$ and $q'$ are quotient maps (by definition),
and since $U_{\bar M} \cup U_{\bar Y}$ 
and $U_{M'} \cup U_{Y'}$ are saturated with respect to the 
maps $\bar q$ and $q'$ (see lemma \ref{lemma:step6}),
$U_{\bar X}$ and $U_{X'}$ are open in $\bar X$ and $X'$.
Therefore $U_x = U_{\bar X} \times U_{X'} \cap X$ is open in $X$
and contains $x$. It is left to show that 
$U_x \subset U$.
Consider $t\in U_x\cap Y$. Then
$q^{-1}t = h^{-1}t \sqcup \{t\} \subset M\sqcup Y$.
But since $j_Y$ is mono,  $t\in (U_{\bar Y} \times U{_Y'})\cap Y$
and therefore $t\in U$.
Otherwise, if $t\in X\minus Y$ then $q^{-1}t = m$, unique 
point belonging to $(U_{\bar M} \times U_{M'}) \cap M \subset q^{-1}U$.
Thus $q^{-1}t \subset q^{-1}U$ and therefore $t\in U$.
We have proved that $U_x \subset U$ and the  proof is now complete.
\end{proof}

\begin{remark}
Since a pull-back of a $\V$-bundle is a $\V$-bundle,
if $M'$ is a $\V$-bundle, then $M$ is a $\V$-bundle.
Furthermore, by 
\ref{lemma:step2}, 
$A$ is the restriction of $M$ to $\bar A$. 
Thus, in this case, proposition \ref{propo:theo:pushout2}
implies that if $\bar M$ and $\bar M'$ are locally finite
countable CW-complexes and  $M'$ a $\V$-bundle, 
if $X'$ is obtained by attaching 
the $\V$-bundle $M'$
to a $\V$-family $Y'$, 
then the pull-back $\V$-family $\bar f^* X'$
is obtained by attaching the pull-back $\V$-bundle $M$ 
to the pull-back $\V$-family $\tilde f^* Y'$.
\end{remark}

\section{Push-out of pull-back $\V$-bundles}
\label{section:pushout}
In this section we complete the proof of the bundle theorem
\ref{theo:bundle}.
In case of the pull-back of $\V$-bundles, 
the following properties hold.

\begin{lemma}\label{lemma:relativeMA}
Let  $\bar Y$ be a space, $\bar h\from S^{n-1} \to \bar Y$ a map
and $\bar X= \bar Y\cup_{\bar h} D^n$. 
Let $\bar \hh\from D^n \to \bar X$ denote the characteristic map.
If $X \to \bar X$ 
is a $\V$-bundle and $Y$ is the restriction
of $X$ to $\bar Y$,
then the following diagram is a push-out
\ldiag{diag:here}{%
\bar h^*Y \ar[d] \ar[r] & 
Y \ar[d] 
\\
\bar \hh^* X \ar[r] & X \\
}%
\end{lemma}
\begin{proof}
Since the diagram commutes, 
one gets the map $t\from P \to X$,
where $P$ denotes the push-out space 
$\bar \hh^*X \cup Y$.
It is easy to see that $t$ is bijective and covers
the identity of $\bar X$.
Now we show that  $t$  is a local homeomorphism
when restricted to 
the space over $\bar U$, 
where $\bar U \subset \bar X$ is an open set such 
that $X|\bar U$ is trivial.
In fact, in this case the diagram \ref{diag:here}
is reduced as 
\ndiag{%
\bar h^* (W\times(\bar U \cap \bar Y)) \ar[d] \ar[r] & 
W \times (\bar U \cap \bar Y)  \ar[d] 
\\
W \times \bar \hh^{-1} \bar U \ar[r] & W \times \bar U, \\
}%
which is a push-out since $W$ is locally compact Hausdorff
and $\bar X = \bar Y \cup_{\bar h} D^n$.
Thus $P|\bar U$ and  $X|\bar U$ are homeomorphic. 
But $P|\bar U$ and $X|\bar U$, when $\bar U$ ranges over 
all the trivializing neighborhoods of $\bar X$, are open
covers of $P$ and $X$ respectively, with 
$t(P|\bar U) = X|\bar U$ for each $\bar U$.
Hence $t$ is an open map.
\end{proof}

Since exactly the same argument can be applied to the attaching
of a set $\bar Z$ of $n$-cells,
we have the following generalization of lemma \ref{lemma:relativeMA}.
\begin{lemma}\label{lemma:attaching}
Let  $\bar Y$ be a space, $\bar Z$ a set 
and $\bar h\from \bar Z \times S^{n-1} \to \bar Y$ an attaching map with
$\bar X= \bar Y\cup_{\bar h} (\bar Z \times D^n)$. 
Let $\bar \hh\from \bar Z \times D^n \to \bar X$ denote the characteristic map.
If $X \to \bar X$ 
is a $\V$-bundle and $Y$ is the restriction
of $X$ to $\bar Y$
then the following diagram is a push-out
\ndiag{%
\bar h^*Y \ar[d] \ar[r] & 
Y \ar[d] 
\\
\bar \hh^* X \ar[r] & X. \\
}%
\end{lemma}

We recall now the following important property.
% \begin{lemma} \label{lemma:limitMn}
Let $M_0 \subset M_1 \subset \dots \subset M_n \subset \lim_{n\geq 0} M_n$
be a sequence of spaces, and $\displaystyle M = \lim_{n\to \infty} M_n$.
Then for every locally compact Hausdorff space $W$
\begin{equation}\label{eq:limitMn}
\lim_{n\to \infty} \left( W \times M_n  \right) =
W \times \left( \lim_{n\to \infty} M_n  \right)  = W \times M.
\end{equation}
%\end{lemma}

\begin{lemma}
\label{lemma:cwpair}
Let  $(\bar M, \bar A)$ be a CW-pair;  let $M \to \bar M$
be a $\V$-bundle and $A$ the restriction of $M$ to
$\bar A$. Then $(M,A)$ is a $\V$-CW-pair.
\end{lemma}
\begin{proof}
Consider the filtration on skeleta
\begin{equation*}
\bar A \subset \bar M_1 \subset \bar M_2 \subset \dots
\subset \bar M_n \subset \dots \subset \bar M.
\end{equation*}
It induces a filtration on restrictions of the bundles
$M\to \bar M$
\begin{equation*}
A \subset  M_1 \subset  M_2 \subset \dots
\subset  M_n \subset \dots \subset  M.
\end{equation*}
By lemma \ref{lemma:attaching} for every $n\geq 0$ 
the diagram
\ndiag{%
\bar h_n^*M_{n-1} \ar[d] \ar[r] & 
M_{n-1} \ar[d] 
\\
\bar \hh_n^* M_n \ar[r] & M_n, \\
}%
is a push-out,
where $\bar h_n\from \bar Z_n \times S^{n-1} \to \bar M_{n-1}$
is the attaching map on $\bar M_{n-1}$ 
and $\bar \hh_n$ is the corresponding characteristic map.
Since $D^n$ is contractible, there exists  a $\V$-isomorphism
of (pair of ) bundles
\begin{equation} \label{eq:isomorph}
(\bar \hh_n^* M_n, \bar h_n^* M_{n-1}) \cong
W \times ( \bar Z_n \times D^n, \bar Z_n \times S^{n-1})
\end{equation}
for a suitable object $W\in \V$.
This means that $M_n$ is obtained by attaching $W \times \bar Z_n \times D^n$
to $M_{n-1}$ with as attaching map the composition of the 
map induced by
pull-back $\bar h^*$ and the isomorphism of \ref{eq:isomorph}.
Thus $(M_n,A)$ there are skeleta $M_n \subset M$
such that 
\begin{equation*}
A \subset  M_0 \subset M_1 \subset \dots M.
\end{equation*}
To show that $(M,A)$ is a relative $\V$-complex (and hence
a $\V$-CW-pair) we need to show that 
$\displaystyle \lim_{n\to \infty} M_n = M$.
Since there is a continuous  bijection $t\from 
\displaystyle \lim_{n\to \infty} M_n \to  M$
covering the identity of $\bar M$, we need to show that 
$M$ has the topology of the limit. It suffices to 
restrict $M$ to the trivializing neighborhoods $\bar U\subset \bar M$,
and this is true as a consequence of \ref{eq:limitMn}.
\end{proof}

Lemma \ref{lemma:cwpair} implies the following corollary,
by taking $A=\emptyset$.
%\begin{coro}
% \label{coro:cwpair}
Let  $\bar M$ be a CW-complex;  let $M \to \bar M$
be a $\V$-bundle. 
Then $M$ is a $\V$-complex.
%\end{coro}
This is the second part of \ref{theo:bundle}.

\section{Push-out of function spaces}
\label{section:9}
We now start the actual proof of the principal bundle
theorem, which is achieved in the next two sections.

Let $V$ be 
an object of $\V$.
We recall that 
a family $\mathcal{K}$ of compact sets of $V$ is termed
\emph{generating} if for every 
$\V$-family $Y$ 
the subsets $N_{K,U}$ with $K\in \mathcal{K}$
and $U$ open in $Y$ yield a sub-basis for the topology of $Y^V$.

\begin{defi}\label{defi:nnep}
Let $Y\subset X$ a pair of $\V$-families and
 a closed inclusion.
We say that $Y$ has the \emph{$N$-neighborhood extension property in $X$}
(NNEP)
if 
for every $V\in \V$
there exists a 
generating family $\mathcal{K}$ of compact subsets of $V$
such that the following is true: 
let $U$ be   open in $X^V$  and 
for  $i=1,\dots, l$ let $K_i \subset V$ be compact sets in $\mathcal{K}$
and $U_i\subset Y$ be open such that 
\begin{equation*}
\bigcap_{i=1}^l N_{K_i,U_i} \subset U \cap Y^V;
\end{equation*} 
then there are $l$ open subsets $U'_i\subset X$ such that 
$U'_i \cap Y = U_i$ and
\begin{equation*}
\bigcap_{i=1}^l N_{K_i,U_i'} \subset U.
\end{equation*}
\end{defi}

An important easy property of the sets $N_{K,U}$ is the following:
if $U\subset Y$  and $h\from A \to Y$ is a map,
then
\begin{equation*}
N_{K,h^{-1}U} = {(h^V)}^{-1} N_{K,U}.
\end{equation*}

\begin{lemma}\label{lemma:coprod}
Let $Z$ be a set and $(M_z,A_z)$ be an NNEP-pair
for every $z\in Z$. Then the coproduct of inclusions
\begin{equation*}
\coprod_{z\in Z} A_z \subset \coprod_{z\in Z} M_z
\end{equation*}
yields a NNEP-pair.
\end{lemma}

\begin{lemma}\label{lemma:nnep}
Assume that $A$ has the NNEP in $M$ and that $h\from A \to Y$
is a $\V$-map to a $\V$-family $Y$. Then 
$Y$ has the NNEP in 
the
push-out space 
$X=M\cup_h Y$. 
\end{lemma}
\begin{proof}
By lemma \ref{lemma:pushout} $Y\to X$ is a closed inclusion.
Let $\hh\from M \to X$ denote the characteristic map of the push-out.
Let $U\subset X^V$  be open, and let $K_i\subset V$ and $U_i \subset Y$
be compact subsets of a generating family $\mathcal{K}$ 
of $V$ and open subsets such that 
\begin{equation*}
\bigcap_{i=1}^l N_{K_i,U_i} \subset U \cap Y^V
\end{equation*}
as in definition \ref{defi:nnep}.
The map $\hh^V$ is continuous, hence $U''={(\hh^V)}^{-1} U \subset M^V$
is open in $M^V$.
Moreover, the intersection
\begin{equation*}
{(h^V)}^{-1} 
\bigcap_{i=1}^l
N_{K_i,U_i} = 
\bigcap_{i=1}^l
N_{K_i,h^{-1}U_i} 
\end{equation*} 
is contained in $U''$.
Since $A$ has the NNEP in $X$, 
there are $l$ open sets $U_i''$ in $M$ such that  
$U_i'' \cap A = h^{-1}U_i$ and 
\begin{equation*}
\bigcap_{i=1}^l N_{K_i,U''_i} \subset U''.
\end{equation*}
Now consider the subsets $U_i'=U_i \bigcup \hh(U_i'') \subset X$. Since
$U_i' \cap Y = U_i$, $\hh^{-1}U'_i = U_i''$
and 
$X$ is the push-out of $M$ and $Y$, 
each $U_i'$ is open in $X$.
Moreover, 
consider 
\begin{equation*}
f\in \bigcap_{i=1}^l N_{K_i,U'_i}.
\end{equation*}
If $f\in Y^V$ then 
\begin{equation*}
f\in \bigcap_{i=1}^l N_{K_i,U'_i} \cap Y^V = 
\bigcap_{i=1}^l N_{K_i,U'_i\cap Y}  = 
\bigcap_{i=1}^l N_{K_i,U_i}  \subset U\cap Y^V, 
\end{equation*}
hence $f\in U$.
On the other hand, if $f\not\in Y^V$, then
there is a unique $f''\from K \to M\minus A$  
such that $f = {\hh^V}(f'')$. Since for every $i$
$fK_i = hf''K_i \subset U'_i$, we have
\begin{equation*}
f''K_i \subset \hh^{-1} U_i' = U''_i,
\end{equation*}
thus
\begin{equation*}
f'' \in \bigcap_{i=1}^l 
N_{K_i,U''_i} \subset U''.
\end{equation*}
Hence $f= \hh^V(f'') \in \hh^V U'' = \hh^V {(\hh^V)}^{-1} U = U$.
Therefore 
\begin{equation*}
\bigcap_{i=1}^l N_{K_i,U'_i} \subset U,
\end{equation*}
and the proof is complete.
\end{proof}

\begin{lemma}
\label{lemma:nnep2}
Assume that $(M,A)$ is a NNEP-pair, $Y$ is a $\V$-family
and $h\from A \to Y$ is a $\V$-map.
Let $X$ be the push-out $X=M\cup_h Y$. Then the 
following diagram
\ndiag{%
A^V \ar[r]^{h^V} \ar@{ >->}[d] 
& 
Y^V \ar@{ >->}[d] \\
M^V \ar[r]^{\hh^V} 
&
X^V
}% 
is a push-out.
\end{lemma}
\begin{proof}
By lemma \ref{lemma:nnep} $(X,Y)$ is a NNEP-pair.
Let $\mathcal{K}$ denote the generating family of compact sets of $V$.
Let $P$ denote the push-out $M^V\cup_{h^V} Y^V$. 
By the push-out property there is a continuous 
map $P \to X^V$. It is easy to see that it is bijective,
so that,
by identifying $P$ and $X^V$, 
the only thing to prove is that if $U$ 
is an open in the push-out topology, then it is open in the compact-open
topology. This is true if and only if 
for every $f_0\in U$ 
there is a subset $U'\subset X^V$ open in the compact-open topology
such that $f_0\subset U' \subset U$.

Consider first the case $f_0\not\in Y^V$. Then there
exists a unique 
$f_0'\in M^V \minus A^V$ such that $\hh^V(f_0') = f_0$.
Since $A^V$ is closed in $M^V$, $f_0'$ is contained in
an open neighborhood $U_0'' \subset M^V$ such that 
$U_0'' \cap A^V = \emptyset$ 
and $U_0''$ is of type
\begin{equation*}
\bigcap_{i=1}^l N_{K_i,U_i''}
\end{equation*}
for some $K_i\subset V$ and some $U_i''$ in $M\minus A$. 
For each $i$ the image $U'_i=\hh U_i''$ is open in $X$,
hence the set
\begin{equation*}
\bigcap_{i=1}^l N_{K_i,U_i'}
\end{equation*}
is open in $X^V$, contains $f_0$ and is contained in $U$.

Now consider $f_0\in Y^V$. The intersection $Y^V\cap U$
is open in $Y^V$ and hence there is a neighborhood of $f_0$
in $Y^V$ contained in $Y^V\cap U$.
Because $Y^V$ has the compact-open topology this means that 
there are $l$ compact subsets $K_i\subset V$ 
and open sets $U_i\subset Y$
such that 
\begin{equation} %\label{eq:defiofN}
f_0 \in
\bigcap_{i=1}^l N_{K_i,U_i} \subset U,
\end{equation}
where as above 
$N_{K_i,U_i}$ denotes the set of all the maps in $Y^V$
such that $fK_i\subset U_i$. 
Without loss of generality
we can assume that $K_i$ belongs to the generating family $\mathcal{K}$
for every $i=1\dots l$.

Consider the sets 
\begin{equation*}
U''=
{(\hh^V)}^{-1} U \subset M^V
\end{equation*}
and 
\begin{equation*}
\bigcap_{i=1}^l N_{K_i,h^{-1}U_i} \subset A^V \cap U''.
\end{equation*}
Since $(M,A)$ is a NNEP-pair,  by definition \ref{defi:nnep}
there are open sets $U_i''\subset M$ 
such that $U_i'' \cap A = h^{-1}U_i$ and 
\begin{equation*}
\bigcap_{i=1}^l N_{K_i,U_i''} \subset U''.
\end{equation*}

Let $U_i'=U_i \bigcup \hh U_i''$ for every $i=1$, $\dots$, $l$.
As in the proof of lemma \ref{lemma:nnep}, 
they are $l$ open subsets of $X$ with the property
that $U_i'\cap Y=U_i$,
hence $f_0\in N_{K_i,U'_i}$ for every $i$.
And, again as in the proof of lemma \ref{lemma:nnep},
\begin{equation*}
U' = \bigcap_{i=1}^l N_{K_i,U_i'} \subset U.
\end{equation*}
This concludes the proof.
\end{proof}

\begin{lemma}\label{lemma:limit}
Consider a sequence of $\V$-families
\begin{equation*}
X_0 \subset X_1 \subset \dots \subset X_n \subset \dots \subset 
X=\lim_{n\to \infty} X_n
\end{equation*}
such that for every $n\geq 1$ $(X_n,X_{n-1})$ is a NNEP-pair.
Then $(X,X_0)$ is a NNEP-pair  and
\begin{equation*}
X^V= \lim_{n\to \infty} X_n^V 
\end{equation*}
\end{lemma}
\begin{proof}
Let $U\subset X^V$ be open,
and $\bigcap_{i=1}^l N_{K_i,U_i^0} \subset U\cap X_0^V$ 
a subset of $X_0^V$,
with $K_i \subset V$ compact containing bases 
and $U_i^0$ open in $X_0$.
By induction, it is possible to define sequences of sets
$U_i^n$ open in $X_n$ such that  for every $\geq 1$
\begin{equation*}
U_i^n \cap X_{n-1} = U_i^{n-1}
\end{equation*}
\begin{equation*}
\bigcap_{i=1}^l N_{K_i,U_i^n} \subset X_n^V \cap U.
\end{equation*}
Now, because $X=
\lim_{n\to \infty}
X_n$,
the sets $U_i^\infty = \bigcup_{n=0}^{\infty} U_i^n$ 
are open in $X$, 
hence 
\begin{equation*}
\bigcap_{i=1}^l N_{K_i,U_i^\infty} \subset U.
\end{equation*}
is an open set. 
Moreover, by construction
$U_i^\infty \cap X_0 = U_i^0$,
therefore $X_0$ has the NNEP in $X$.
The same construction shows that $X^V = \lim_{n\to \infty} X_n^V$.
\end{proof}

\section{NKC-categories}

\label{section:10}
The crucial property of a 
NKC-category as defined in 
\ref{defi:NKC}
is the next result.

\begin{propo}
\label{propo:theo:tech} 
Assume that $\V$ is NKC.
Then 
for every $W\in \V$ the pair $(W\times D^n, W \times S^{n-1})$
is a NNEP-pair.
\end{propo}

The proof of proposition \ref{propo:theo:tech} will take the rest of the section.

% from here set junk
\begin{lemma}\label{lemma:product}
Let $R$ be a space and $V$ and $W$ objects of $\V$.
Then 
$$(W\times R)^V = W^V \times R.$$
\end{lemma}
\begin{proof}
Let $\pr_1$ and $\pr_2$ denote the projections onto the 
first and the second factor of $W\times R$.
The map $\pr_2\epsilon_0\from (W\times R)^V \to R$ 
sending $f\from V \to W\times R$
to $\pr_2 f(0) \in R$ is continuous. 
Furthermore, the map $\pr_1^V\from (W\times R)^V \to W^V$
is continuous.
Therefore the map $F=(\pr_1^V, \pr_2\epsilon_0)\from 
(W\times R)^V \to W^V \times R$ is continuous.
Consider now the function $G\from W^V \times R \to (W\times R)^V$
defined by $G(f,r)(v) = (f(v),r)$ for every $v\in V$ and 
every $r\in R$. It is the adjoint of 
the map $V\times W^V \times R \to W \times R$ 
defined by $(v,f,r) \mapsto (f(v),r)$, which is continuous
since the evaluation map is continuous. Therefore $G$ 
is continuous.
It is readily seen that $GF=1$ and $FG=1$.
\end{proof}

Let $D^n_\epsilon = \{ x\in D^n \st |x|>1-\epsilon \} $,
for any $\epsilon$ small,
and let $\rho_{\epsilon}\from  D^n_\epsilon  \to S^{n-1}$ 
be the retraction defined by $x \mapsto \frac{x}{|x|}$.
With an abuse of notation we use the same symbol for the 
induced retraction
$\rho_\epsilon=1_{W} \times 
\rho_\epsilon\from W\times D_\epsilon^n \to W\times S^{n-1}$.
For every $W\in \V$ there is an induced retraction
\begin{equation*}
\rho_\epsilon^V\from (W\times D^n_\epsilon) ^V
\to
(W\times S^{n-1})^V.
\end{equation*}

\begin{lemma}\label{lemma:epsilonk}
If $C$ is a compact in $(W\times S^{n-1})^V$,
and $C$ is contained in an open set $A$ of $(W\times D^n)^V$,
then there exists $\epsilon>0$ such that 
${(\rho^{V}_\epsilon)}^{-1}C \subset A$.
\end{lemma}
\begin{proof}
If $\epsilon$ is small, we have $D^n_\epsilon = S^{n-1}\times 
[0,\epsilon)$, and hence 
$ (W\times D^n_\epsilon) ^V = (W\times S^{n-1})^V \times [0,\epsilon)$
and $\rho^V_\epsilon$ is the projection onto the first factor.
Each $x\in C$ is a point of $A$, therefore there exists an 
open neighborhood of $x$ 
of type $O_x\times [0,\epsilon_x)$ contained in $A$.
Being $C$ compact there is a finite number of 
points $x$ such that 
$C$ is covered by $O_x\times [0,\epsilon_x)$. Therefore
there exists a $\epsilon>0$ such that 
$\cup_{x} O_x \times [0,\epsilon)$ is an open
subset of $A$ containing $C$ and   
thus ${(\rho^V_\epsilon)}^{-1} C \subset A$.
\end{proof}

Now we can start the proof of  
proposition \ref{propo:theo:tech}. 
Let $U$ be an open subset 
$U\subset (W\times D^n)^V$.
Let $\mathcal{K}$ be the generating family of compact
subsets of $V$ of  the NKC property.
Let $K_i\subset V$ be compact subsets  in $\mathcal{K}$
and $U_i\subset W\times S^{n-1}$
open subsets,
such that 
\begin{equation}\label{eq:defiofN}
\bigcap_{i=1}^l N_{K_i,U_i} \subset U \cap (W\times S^{n-1})^V.
\end{equation}
By lemma \ref{lemma:a1} there is a 
sequence of open sets $A^k_i \subset W\times S^{n-1}$
with $k\geq 1$ such that, for every $k$,
$A_i^k \subset A^{k+1}_i$,
\begin{equation}\label{eq:this}
\bigcup_{k\geq 1} A_i^k = U_i
\end{equation}
and the closure $\bar A_i^k$ is a compact subset of $U_i$.

\begin{lemma}\label{lemma:Ncompact}
For every $i=1,\dots, l$ the space
$N_{K_i,\bar A_i^k}$ is compact.
\end{lemma}
\begin{proof}
It is closed, because 
it is equal to the intersection
of all the spaces $N_{\{x\},\bar A_i^k}$, with $x\in K_i$,
and they are closed because pre-images of the closed set
$\bar A_i^k$ under the 
evaluation  map  
$(W\times S^{n-1})^V\times \{x\} \to W\times S^{n-1}$,
which is continuous because of lemma \ref{coro:adjoint}.
Now we prove that it is compact.
Consider the projection $p\from W\times S^{n-1} \to W$ 
onto the first factor. Since $\bar A_i^k$ is compact, its 
image $p\bar A_i^k$ is compact in $W$.
Moreover, it is easy to see that 
\[
p^V(N_{K_i,\bar A_i^k}) \subset  
N_{K_i,p\bar A_i^k},
\]
and hence that 
\[
N_{K_i,\bar A_i^k} \subset 
{(p^V)}^{-1}N_{K_i,p\bar A_i^k} = N_{K_i,p\bar A_i^k} \times S^{n-1}.
\]
Now, since $K_i$ and $p\bar A_i^k$ are compact
and $\V$ is an NKC-category,
$N_{K_i,p\bar A_i^k}$ is compact. 
Therefore $N_{K_i,\bar A_i^k}$ is a closed subset 
of a compact space, and hence compact. 
\end{proof}

Consider for every $k\geq 1$ the set
\begin{equation}
C_k = \bigcap_{i=1}^l N_{K_i,A^k_i} \subset (W\times S^{n-1})^V.
\end{equation}
It is open, by definition of compact-open topology.
Let $C_{\infty}$ denote the union $\bigcup_{k\geq 1} C_k $.

\begin{lemma}\label{lemma:a3}
We have that 
\begin{equation*}
C_\infty = 
\bigcap_{i=1}^l N_{K_i,U_i} \subset (W\times S^{n-1})^V \cap U.
\end{equation*}
\end{lemma}
\begin{proof}
A $\V$-map $f\in (W \times S^{n-1})^V$ belongs to $C_{\infty}$
if and only if there is $k\geq 1$ such that 
$f\in C_k$, that is if and only if there is $k\geq 1$ 
such that $fK_i \subset A^k_i$ for every $i=1,\dots, l$.
Of course this implies that $fK_i \subset U_i$ for every $i$.
On the other hand, if $fK_i \subset  U_i$, then 
by equation \ref{eq:this} there is $k_i\geq 1$ such  that 
$fK_i \subset A^k_i$ for every $k\geq k_i$. This implies
that there is $k$ (take the maximum of all $k_i$) 
such that $fK_i \subset A^k_i$ for 
every $i=1,\dots, l$.
\end{proof}

\begin{lemma}\label{lemma:closure}
The closure $\bar C_k \subset (W\times S^{n-1})^{V}$ 
is compact and contained in $C_{\infty}$.
\end{lemma}
\begin{proof}
Because  of its 
definition $\bar C_k \subset \cap_{i=1}^l \overline{N_{K_i,A_i^k}}$.
Because for every $i=1,\dots, l $ the closure
$\overline{N_{K_i,A_i^k}}$ is contained in 
$N_{K_i,\bar A^k_i}$, which is compact 
by lemma
\ref{lemma:Ncompact},
we have that $\bar C_k$ is compact.
Furthermore, by assumption $\bar A_i^k$ is contained in $U_i$,
and hence $N_{K_i,\bar A_i^k}$
is contained in 
$N_{K_i, U_i}$.  Therefore the closure of $C_k$  
is contained in $C_{\infty}$.
\end{proof}

Because of lemma \ref{lemma:closure} and lemma \ref{lemma:a3}
the compact set 
$\bar C_k$  is contained  in 
$U$ 
in $(W\times D^{n})^V$.
By lemma \ref{lemma:epsilonk} it is possible 
to find an $\epsilon_k>0$ 
such that 
\begin{equation}\label{eq:contained2}
{(\rho_{\epsilon_k}^V)}^{-1} C_k \subset 
U.
\end{equation}
Without loss of generality we can assume that the sequence
$\epsilon_k$, with $k\geq 1$, is decreasing.

Now consider for $i=1,\dots, l$ the following 
sets in $(W\times D^n)^V$:
\begin{equation}\label{eq:defi}
U_i' =  
\bigcup_{k\geq 1} 
\rho_{\epsilon_k}^{-1} A^k_i.
\end{equation}
Since it is the union of open sets, 
for every $i=1,\dots, l$, the set $U'_i$  is open in $W\times D^n$.
Hence the set
\begin{equation}
U' = \bigcap_{i=1}^l N_{K_i,U_i'} \subset (W\times D^n)^V,
\end{equation}
is open.

\begin{lemma}\label{lemmaeq:4}
For every $i=1,\dots, l$,
$%\begin{equation*}
U_i' \cap (W\times S^{n-1}) = U_i
$.%\end{equation*}
\end{lemma}
\begin{proof}
Because of \ref{eq:defi}, 
it suffices to prove that  for every $i$
and every $k\geq 1$
\begin{equation*}
\rho_{\epsilon_k}^{-1} A^k_i \cap (W\times S^{n-1}) \subset U_i,
\end{equation*}
and this is true because $\rho_{\epsilon_k}^{-1} A^k_i \cap (W\times S^{n-1})
= A^k_i$.
\end{proof}

\begin{lemma}\label{lemma:contained}
$U'\subset U$.
\end{lemma}
\begin{proof}
If $f\in U'\cap (W\times S^{n-1})^V$, 
then $fK_i \subset (W\times S^{n-1}) \cap U'_i$,
and by lemma \ref{lemmaeq:4} this implies $fK_i \subset U_i$. Thus,
by \ref{eq:defiofN}, 
$f_0 \in  U$.
On the other hand, assume that 
$f\in (W\times e^n)^V$. 
For every $i=1,\dots, l$ 
the  image of the compact $f(K_i)$ is in $U_i'$ and hence, 
because of \ref{eq:defi}, for every $i$ 
there exists an integer  $k_i\geq 1$ (depending on $i$) such  
$f'(K_i) \subset \rho_{\epsilon_{k_i}}^{-1} A_i^{k_i}$. 
In particular, for every $i=1,\dots, l$,
\begin{equation*}
f'(K_i) \subset W\times D^n_{\epsilon_{k_i}},
\end{equation*}
and thus, because $f$ is a $\V$-map and $\epsilon_k$ 
is a decreasing sequence,
\begin{equation*}
f'(V) \subset W\times D^n_{\epsilon_{m}},
\end{equation*}
where $m$ denotes the maximum of $k_1,\dots, k_l$.
Moreover, the sequence of spaces $A_i^k$
is increasing with $k\geq 1$, and therefore 
$A_i^{k_i} \subset A_i^m$ for every $i=1,\dots, l$.
Thus for every $i=1,\dots, l$,
\begin{equation*}
f'(K_i) \subset \rho_{\epsilon_{m}}^{-1} A_i^{m}. 
\end{equation*}
But this means that for every $i$
\begin{equation*}
f \in N_{K_i,\rho_{\epsilon_m}^{-1}A_i^m} =
{(\rho_{\epsilon_m}^V )}^{-1} N_{K_i,A_i^m},
\end{equation*}
hence that 
\begin{equation*}
f \in (\rho_{\epsilon_m}^V )^{-1} C_m.
\end{equation*}
Because of equation \ref{eq:contained2} this implies that 
$f\in U$. 
As claimed, we have proved that $U' \subset U$.
\end{proof}

This is the end of the proof 
of proposition \ref{propo:theo:tech}.
We can now draw the consequences which are needed.

%%%%
%%%%
%%%%
%%%%
%%%%
%%%%

\begin{lemma}\label{lemma:step}
Let $\V$ be a NKC-category.
Let $Y$ be a $\V$-family, $Z$ a $\V$-set 
and $h\from Z\times S^{n-1} \to Y$
a $\V$-map, with $n\geq 1$. 
If $X$ is a family obtained as a push-out
\ndiag{%
Z \times S^{n-1} \ar[r]^{\hspace{24pt} h} \ar@{ >->}[d] 
& Y \ar@{>->}[d] \\ 
Z\times D^n \ar[r]^{\hspace{12pt} \hh} 
& X \\
}%
then  the following diagram
\ndiag{%
Z^V\times S^{n-1} \ar@{=}[r] \ar@{ >->}[d]&
(Z\times S^{n-1})^V \ar[r]^{\hspace{24pt} h^V} \ar@{ >->}[d] 
& Y^V \ar[d] \\
Z^V \times D^n \ar@{=}[r] & 
(Z\times D^n)^V \ar[r]^{\hspace{24pt} \hh^V} & X^V\\
}%
is a push-out,  
and $Y$ has the $N$-neighborhood extension property in $X$.
\end{lemma}
\begin{proof}
By lemma \ref{lemma:coprod} and proposition \ref{propo:theo:tech}
if $Z$ is a $\V$-set the pair $(Z\times D^n, Z \times  S^{n-1})$
is a NNEP-pair.
Hence, by lemma \ref{lemma:nnep} and lemma \ref{lemma:nnep2}
$X^V$ is the push-out of $(Z\times D^n)^V$ and $Y^V$
and $Y$ has the $N$-neighborhood extension property in $X$.
\end{proof}

\begin{coro}\label{coro:main}
Let $(X,D)$ be a relative $\V$-complex 
where $\V$ is a NKC-category.
Let $\hh: Z_n\times D^n \to X$
be its characteristic map of $n$-cells and let  $h_n$ denote
the $n$-attaching map of $X$, i.e.~the restriction of $\hh$
to $Z_n\times S^{n-1}$.
For every $n\geq 1$ and every $V\in \V$ the diagram
\ndiag{%
Z_n^V\times S^{n-1} \ar@{=}[r] \ar@{ >->}[d]&
(Z_n\times S^{n-1})^V \ar[r]^{\hspace{24pt} h_n^V} \ar@{ >->}[d]
& X_{n-1}^V \ar[d] \\
Z_n^V \times D^n \ar@{=}[r] &
(Z_n\times D^n)^V \ar[r]^{\hspace{24pt} \hh^V} & X_n^V\\
}%
is a push-out (in $\top$). 
Moreover, 
\begin{equation*}
X^V =
 \lim_{n\to \infty} X_n^V. 
\end{equation*}
\end{coro}
\begin{proof}
The first part 
is a direct consequence of lemma \ref{lemma:step}.
Then, by applying \ref{lemma:limit} we obtain that 
$(X,X_n)$ is a NNEP-pair for every $n$ 
and that $X^V= \lim_{n\to \infty} X^V_n$.
\end{proof}

\begin{coro}
\label{coro:relat}
Let   $(X,D)$ be a relative $\V$-complex (or 
a $\V$-CW-pair) with $\V$ NKC. 
Then $(X,D)$ is a NNEP-pair.
\end{coro}

\begin{coro}\label{coro:nkc} %\label{coro:cwpair}
Let $\V$ be NKC.
Let $(M,A)$ be a $\V$-CW-pair, $Y$ a $\V$-family and $h\from A \to Y$
a $\V$-map. If $X$ is the push-out $X=M\cup_h Y$  then
$(X,Y)$ is a NNEP-pair.
\end{coro}
\begin{proof}
By corollary \ref{coro:relat} $(M,A)$ is a NNEP-pair,
hence by lemma \ref{lemma:nnep} and lemma \ref{lemma:nnep2}
$(X,Y)$ is a NNEP-pair.
\end{proof}

\section{Examples of NKC categories}
\label{section:NKC}

\begin{propo}\label{propo:theo:allcompact}
If 
the fibre functor $\ff$ is faithful and the 
Hom-sets $\hom_\V(V,W)$  
are compact spaces
then $\V$ is a NKC category.
\end{propo}
%\begin{proof}
%Since fibres in $\V$ are by definition Hausdorff, 
%for every $V$, $W\in \V$, every compact 
%$K\subset V$ and compact $C\subset W$ the 
%set $N_{K,C}$ is a closed subset of $\ff \hom_\V(V,W)$.
%Hence if $\ff \hom_\V(V,W)$ is a compact space for all $V$, $W$,
%then $\V$ is a NKC-category. 
%\end{proof}

The proof of the previous proposition is simple. 
The interesting fact is that there are NKC structure categories
with non-compact hom-sets, as a consequence of the following 
proposition.

\begin{propo}\label{theo:NKC} % ex: propo:theo:NKC
If $\V$ is a closed subcategory of the category $\vect$
of finitely dimensional $\R$-vector spaces and linear maps,
then $\V$ has the NKC property.
\end{propo}
\begin{proof}
See \cite{sm}, lemma  (4.2).%\ref{sm_lemma:NKC}.
\end{proof}

This result yields the principal stratified bundle
associated to a $\V$-stratified vector bundle.

%======================================================================

\end{document}